\documentclass[11pt,a4paper]{article}
\usepackage[T1]{fontenc}
\usepackage{lmodern}
\usepackage[utf8]{inputenc}
\usepackage[margin=24mm,top=23mm,bottom=24mm]{geometry}
\usepackage{amsmath,amssymb,amsthm,mathtools}
\usepackage{microtype}
\usepackage[hidelinks,pdfencoding=auto]{hyperref}
\hypersetup{pdfauthor={Jin Hong and Yong Wang},pdftitle={Trilinear noncommutative residues on odd-dimensional manifolds with boundary}}
\newtheorem{theorem}{Theorem}[section]

\newtheorem{proposition}[theorem]{Proposition}
\newtheorem{corollary}[theorem]{Corollary}
\theoremstyle{remark}

\numberwithin{equation}{section}
\allowdisplaybreaks[2]
\title{Trilinear noncommutative residues on odd-dimensional manifolds with boundary}
\author{%
  Jin Hong\textsuperscript{1,2}\qquad
  Yong Wang\textsuperscript{2,}\thanks{Corresponding author.\newline
    \textit{Email addresses:}
    \href{mailto:jhong@nenu.edu.cn}{\nolinkurl{jhong@nenu.edu.cn}} (Jin Hong),
    \href{mailto:wangy581@nenu.edu.cn}{\nolinkurl{wangy581@nenu.edu.cn}} (Yong Wang).}\\[0.6em]
  \parbox{0.94\textwidth}{\centering\small
    \textsuperscript{1}School of Mathematics and Statistics, Institute of Applied Mathematics,\\
    Yili Normal University, Yining 835000, China\\[0.25em]
    \textsuperscript{2}School of Mathematics and Statistics, Northeast Normal University, Changchun, 130024, China%
  }%
}
\date{}
\begin{document}
\maketitle
\vspace{-2.2em}
\begin{abstract}
We compute two factorizations of a trilinear noncommutative residue on a
compact Riemannian manifold of dimension $2m+1\geq5$ with boundary. The
operators are a Dirac operator with torsion, a Connes-type operator on
twisted spinors, and the de Rham--Hodge operator on differential forms.
The interior residue density vanishes in each case. For a warped collar,
both boundary densities contain the same normal derivative of a metric
contraction of the three vector fields, together with the corresponding
mean-curvature term. Their coefficients have ratio $-(2m-1)/3$.
The two factorizations give the same coefficient for the alternating
three-form contribution. A weighted sum therefore recovers the skew
torsion for the torsion Dirac operator and a traced three-form for the
Connes-type operator. In odd dimension the latter also receives a
contribution from the Hodge dual of the perturbing form of degree $2m-2$.
For the de Rham--Hodge operator the weighted sum vanishes, and a normal
insertion in either residue recovers the total mean curvature up to an
explicit constant.
\end{abstract}

\section{Introduction}\label{sec:intro}

Spectral geometry asks how much of a geometric background can be read from the operators naturally attached to it. In this question the noncommutative residue has a particularly local character: it extracts a distinguished homogeneous symbol coefficient and converts it into an invariant integral. The gravitational interpretation of residues of inverse Dirac powers, exemplified by the work of Kalau and Walze \cite{KW}, shows that this construction can recover familiar geometric action densities. The spectral action principle \cite{CC} places the relation between Dirac operators and physical actions in a broader setting. These developments motivate a more refined question: can operator residues distinguish the different tensorial fields that enter a Dirac-type operator, rather than only produce a scalar invariant of their combined geometry?

Clifford insertions provide a natural way to formulate this question. A scalar residue integrates over all directions, whereas multiplication by prescribed vector fields retains directional information before the fiber trace is taken. Three Clifford factors are especially suited to an alternating three-form, because their product contains a degree-three component together with metric contractions of degree one. The coexistence of these components is essential here. It makes a trilinear residue sensitive both to skew torsion and to normal geometric data, and it suggests that suitable combinations of residues may separate the two. The inserted fields are therefore probes of the geometry; they are not assumed to satisfy dynamical equations.

The boundary changes this problem in a fundamental way. The geometry now includes the embedding of the boundary, the normal variation of the metric, and the way the probing fields extend into a collar. Moreover, truncation does not commute with multiplication of pseudodifferential operators. As a result, the algebraic factorization of an inverse power, immaterial in the interior modulo smoothing operators, becomes part of the boundary construction. It is precisely this dependence that we exploit. Instead of viewing the difference between factorizations as an ambiguity to be removed, we compare two fixed factorizations and use their distinct coefficients to distinguish geometric contributions.

There is a complementary physical reason to retain boundary information. The gravitational action on a region with boundary involves extrinsic curvature, as in the work of York and of Gibbons and Hawking \cite{York,GH}. Boundary terms also arise naturally in the spectral action \cite{CCB,CCBoundary}. Heat-kernel methods provide a systematic account of local geometric contributions in effective-action calculations \cite{Vass}. These connections motivate our interpretation, but the analytic objects must be distinguished: the residues considered here are defined by truncated parametrices in the boundary calculus, rather than by the spectrum of a specified self-adjoint boundary-value problem. Recovering the same integrated mean-curvature invariant does not, by itself, identify a complete gravitational action or a boundary condition for a quantum field.

Skew torsion supplies the second geometric ingredient. Metric connections with totally antisymmetric torsion occur in non-integrable geometric structures and in spinorial constructions related to mathematical physics \cite{Agricola}. Their associated Dirac operators contain a zero-order Clifford three-form, making them a natural testing ground for tensor recovery. The spectral treatment of gravity with torsion \cite{PS} and the trilinear spectral torsion functional \cite{DSZ} suggest that residues can distinguish information that is not captured by the Levi--Civita metric alone. Our aim is to understand this distinction at a boundary, where the three-form coupling competes with first normal derivatives and extrinsic curvature.

A trilinear insertion of Clifford multiplications can recover a tensor from
the residue of a Dirac operator. For a closed spin manifold, the spectral
torsion functional of D\k{a}browski, Sitarz and Zalecki recovers the
totally antisymmetric torsion of a metric connection \cite{DSZ}. On a
manifold with boundary the choice of how to factor the inverse power
matters. Truncation to the manifold and composition produce a singular
Green term in Boutet de Monvel's calculus \cite{BdM}. The trace constructed
by Fedosov, Golse, Leichtnam and Schrohe assigns a local boundary density
to this term \cite{FGLS}. Two products with the same interior symbol can
thus yield different geometric information at the boundary.

Let $M$ have dimension $n=2m+1$, $m\geq2$, and let $X,Y,Z$ be smooth
vector fields on $M$. For each operator $\mathcal D$ considered below,
write $A=c(X)c(Y)c(Z)$ on its Clifford module. We study
\begin{gather}
 \widetilde{\operatorname{Wres}}\!\left[
       \pi^+(A\mathcal D^{-1})\circ\pi^+(\mathcal D^{-2m+2})\right],
 \label{eq:def-first}\\
 \widetilde{\operatorname{Wres}}\!\left[
       \pi^+(A\mathcal D^{-2})\circ\pi^+(\mathcal D^{-2m+3})\right].
 \label{eq:def-second}
\end{gather}
The operator in the interior of both products is $A\mathcal D^{-(2m-1)}$.
Its symbol of degree $-(2m+1)$ is odd in the cotangent variable. This
parity makes the boundary term the entire residue and allows the effect
of changing the factorization to be read directly from the boundary
density.

Wang obtained boundary residue formulas for the Dirac and signature
operators and related them to the gravitational boundary term
\cite{Wang}. Wang, Wang and Yang treated Dirac operators with torsion
\cite{WWY}. Yang and Wang studied trilinear expressions involving
one-forms and the ordinary Dirac operator in odd dimension \cite{YW}.
In even ambient dimension, Wang and Wang introduced the Connes-type
operator on a twisted spinor bundle and computed the boundary residue
for a perturbation by a three-form and a one-form
\cite[Definition~2.2 and Theorem~4.14]{WW}. Their first splitting uses
the powers in \eqref{eq:def-first}. Passing from dimension $2m$ to
$2m+1$ lowers the required boundary symbol by one further order. This
brings in the normal derivatives of the insertion and the subleading
symbols of both inverse powers.

There are two distinct contributions. The metric contribution depends on
the full first normal derivative of
$g(X,Y)Z_n-g(X,Z)Y_n+g(Y,Z)X_n$ and on the trace of the second
fundamental form. The three-form contribution depends only on boundary
values. For example, for the torsion Dirac operator
$D_T=D+\tfrac32c(T)-mc(V)$, the formulas in Section~\ref{sec:torsion}
give the identity
\begin{equation}
\begin{aligned}
 &\frac{3}{2m+2}\widetilde{\operatorname{Wres}}\!\left[\pi^+\!\left(AD_T^{-1}\right)\circ\pi^+\!\left(D_T^{-2m+2}\right)\right]\\
 &\quad{}+\frac{2m-1}{2m+2}\widetilde{\operatorname{Wres}}\!\left[\pi^+\!\left(AD_T^{-2}\right)\circ\pi^+\!\left(D_T^{-2m+3}\right)\right]\\
 &\qquad=-\frac{6\pi(2m-2)!\operatorname{Vol}(S^{2m-1})}
                {2^m m!(m-2)!}
         \int_{\partial M} T(X,Y,Z)\,d\operatorname{vol}_{\partial M}.
\end{aligned}
\label{eq:torsion-main}
\end{equation}
Thus the combination on the left is tensorial in the boundary values
of the insertions, even though the two separate residues involve their
normal derivatives. Localizing one insertion by a smooth function gives
the boundary values of $T$ from these integrals. Allowing an insertion
with a normal component also detects $T(\partial_{x_n},Y,Z)$.
For the Connes-type operator the corresponding tensor is a trace-weighted
combination of two complementary form degrees. For differential forms
with the Levi--Civita connection the metric contribution is the only one.

The comparison has three consequences that clarify the scope of the results. First, the cancellation of the metric contribution holds for arbitrary smooth extensions of the probing fields. The resulting weighted functional is therefore tensorial in their boundary values, although neither residue separately has this property. Second, the complementary difference removes the alternating potential. With a controlled normal insertion, it recovers total mean curvature independently of the skew torsion. Third, replacing the spinor module by the full exterior module preserves the geometric structure of the answer but changes its normalization through the fiber dimension. This last comparison exhibits a boundary invariant available without a spin structure.

Odd dimension introduces a further representation-theoretic feature for the Connes-type operator. The Clifford volume element acts as a scalar on the chosen irreducible spinor module. Consequently, complementary form degrees can have the same Clifford action up to a fixed phase. The alternating tensor selected by our residues is therefore the trace-weighted combination of degrees three and $n-3$, not the degree-three component in isolation. This is an explicit limitation on reconstruction as well as a useful description of what the observable actually measures. In particular, invisibility of a potential to these two residues is weaker than vanishing of that potential as an operator.

The calculations are carried out for the warped collar in \eqref{eq:collar}. This hypothesis gives a controlled family of hypersurfaces and makes their extrinsic geometry depend on a single scalar normal derivative. Our conclusions about mean curvature and half-density variation are made within this setting; the formulas do not assert that a general, non-umbilic boundary has no additional contributions. Similarly, the parity argument removes the particular interior residue density studied here, not the interior geometry of the operator. These distinctions allow the geometric interpretation to remain closely tied to the actual symbol calculation.

Section~\ref{sec:prelim} fixes the operator, trace and collar conventions.
Section~\ref{sec:interior} proves interior vanishing. Sections
\ref{sec:torsion}--\ref{sec:hodge} compute both boundary densities
for the three operators and identify the mean-curvature and torsion
terms that can be recovered from them.

\section{Preliminaries}\label{sec:prelim}

\subsection{The collar and the three Clifford modules}
Let $(M,g)$ be compact and oriented, with smooth boundary. Assume
that $M$ is spin whenever spinors are used. On a collar of the
boundary let $x_n$ be inward normal distance, with
\begin{equation}
 g=dx_n^2+h(x_n)^{-1}g^{\partial M},\qquad h(0)=1,\qquad h>0 .
 \label{eq:collar}
\end{equation}
At $x_0\in\partial M$ take tangential geodesic coordinates and an
orthonormal boundary frame parallel at $x_0$ for the induced
connection. Extend the frame by parallel transport along the normal
geodesics and set $e_n=\partial_{x_n}$. Write $X_a=g(X,e_a)$ for
$1\leq a\leq 2m$ and $X_n=g(X,\partial_{x_n})$, and use the same
convention for $Y,Z,V$. If $\xi=\xi'+\xi_n dx_n$, then at $x_0$
\[
 c(\xi')=\sum_{a=1}^{2m}\xi_a c(e_a),\qquad
 c(U)c(W)+c(W)c(U)=-2g(U,W).
\]
Here $c(U)$ on a vector is Clifford multiplication by its metric dual.
With the inward convention the unnormalized mean curvature is
\begin{equation}
 K=\sum_{a=1}^{2m}g(\nabla^L_{e_a}\partial_{x_n},e_a)
   =-m h'(0).
 \label{eq:curvature}
\end{equation}
Indeed, the second fundamental form is one half of the normal derivative
of the tangential metric in \eqref{eq:collar}.

On the complex spinor bundle $S(TM)$ let
$D=\sum_{j=1}^n c(e_j)\nabla^{S(TM)}_{e_j}$ be the Levi--Civita
Dirac operator. We use the following torsion convention \cite{WWY,WW,PS}.
\begin{equation}
 D_T=D+\frac32c(T)-m c(V),\qquad
 c(T)=\sum_{1\le j<k<l\le n}
 T(e_j,e_k,e_l)c(e_j)c(e_k)c(e_l).
 \label{eq:torsion-def}
\end{equation}
Here $T$ is a real three-form and $V$ is the vector field describing
vectorial torsion. Formula \eqref{eq:torsion-def} fixes the normalization
of both components. The integer powers used below are defined by
elliptic parametrices, including when vectorial torsion is present.

Let $F\to M$ be a Hermitian vector bundle with a unitary connection.
With the tensor product connection on $S(TM)\otimes F$, define
\begin{equation}
 D_F=\sum_{j=1}^{n}c(e_j)\nabla^{S(TM)\otimes F}_{e_j},
 \qquad \widetilde D_F=D_F+c(\omega)\otimes\Phi,
 \qquad \Phi\in\Gamma(\operatorname{End}F),\quad\Phi^*=\Phi .
 \label{eq:connes-def}
\end{equation}
We use the differential expression of \cite[Definition~2.2]{WW}
on the odd-dimensional spinor module. A smooth complex form is
written $\omega=\sum_{q=0}^n\omega_{[q]}$. Its Clifford action is
\[
 c(\omega)=\omega_{[0]}+
 \sum_{q=1}^n\ \sum_{1\le j_1<\cdots<j_q\le n}
 \omega_{[q]}(e_{j_1},\ldots,e_{j_q})c(e_{j_1})\cdots c(e_{j_q}).
\]
The metric and orientation define the complex-linear Hodge star $*$ on
ambient forms. In a positively oriented orthonormal frame our odd
spinor convention is
$i^{m+1}c(e_1)\cdots c(e_n)=1$.
Write $\kappa_n=i^{-(m+1)}$ for the scalar by which this Clifford
product acts. The rank $\operatorname{rk}(F)=\dim_{\mathbb C}F_x$
is locally constant on $M$. The symbol $\operatorname{tr}_F(\Phi)$ denotes
the ordinary trace of the endomorphism at $x$. One may work on each
connected component separately.

Finally, let $D_{\mathrm{dR}}=d+\delta$ on the full complex exterior
bundle $\Lambda^*T^*M\otimes\mathbb C$, where $\delta$ is the formal
adjoint of $d$ \cite[Sec.~2]{LWW}. Exterior multiplication by $g(U,\cdot)$ minus contraction
with $U$ gives the action $c(U)$, and their sum gives
$\widehat c(U)$. The Levi--Civita connection induced on forms gives
$D_{\mathrm{dR}}=\sum_{j=1}^n c(e_j)\nabla^{\Lambda^*T^*M}_{e_j}$.
The trace is the ordinary trace on the full module in each case.
The respective fiber dimensions are $2^m$, $2^m\operatorname{rk}(F)$,
and $2^{2m+1}$. Scalar symbols act by scalar multiplication, and on
twisted spinors $A$ acts through the spinor factor.

\subsection{The local boundary formula}
Extend the coefficients smoothly across the boundary, and take negative
powers on a closed extension using a parametrix modulo smoothing
operators. The notation $\pi^+P$ denotes the truncation of $P$ to $M$.
On rational symbols in the normal covariable, $\pi^+_{\xi_n}$ retains
the principal parts at poles in the upper half-plane. The symbol product
uses $D_x^\alpha=(-i)^{|\alpha|}\partial_x^\alpha$ and
\[
 \sigma(PQ)\sim\sum_{\alpha}
 \frac{1}{\alpha!}\,\partial_\xi^\alpha\sigma(P)\,
                      D_x^\alpha\sigma(Q).
\]
We fix the scale of $\widetilde{\operatorname{Wres}}$ by taking
coefficient $1$ in front of the integral of the interior symbol on the
unit cosphere. In the boundary formula for a product of truncations
\cite{FGLS,Wang}, the contribution indexed by $(r,\ell,j,k,\alpha)$
to the boundary scalar is
\begin{equation}
\begin{aligned}
 &\frac{(-i)^{|\alpha|+j+k+1}}{\alpha!(j+k+1)!}
 \int_{|\xi'|=1}\int_{\mathbb R}
 \operatorname{tr}_E\Big[
   \partial_{x_n}^j\partial_{\xi'}^\alpha\partial_{\xi_n}^k
                      \pi^+_{\xi_n}\sigma_r(P)\\[-1pt]
 &\hspace{28mm}{}
   \cdot\partial_{x'}^\alpha\partial_{\xi_n}^{j+1}\partial_{x_n}^k
                      \sigma_\ell(Q)\Big](x_0,\xi)\,d\xi_n\,dS(\xi')
 .
\end{aligned}
\label{eq:boundary-formula}
\end{equation}
Here $r+\ell-j-k-|\alpha|-1=-n$. The integers $j,k$ are nonnegative,
$\alpha$ is a tangential multi-index, and $r,\ell$
are bounded above by the respective orders. The sum over admissible
tuples, multiplied by $d\operatorname{vol}_{\partial M}$, is the
boundary density. Its integral over $\partial M$ is the boundary
contribution to the residue. For \eqref{eq:def-first} and
\eqref{eq:def-second} precisely five cases meet this order condition.
In cases $a_1,a_2,a_3$ both factors contribute their principal symbols,
so $(r,\ell)=(-1,-2m+2)$ in \eqref{eq:def-first} and
$(r,\ell)=(-2,-2m+3)$ in \eqref{eq:def-second}.
Case $a_1$ has $j=k=0$ and $|\alpha|=1$, with coefficient $-1$.
Case $a_2$ has $j=1$ and $k=|\alpha|=0$, while case $a_3$ has
$k=1$ and $j=|\alpha|=0$. Each has coefficient $-1/2$.

In cases $b$ and $c$, $j=k=|\alpha|=0$ and the coefficient is $-i$.
For \eqref{eq:def-first}, the right symbol drops by one order in case $b$,
giving $(r,\ell)=(-1,-2m+1)$, and the left symbol drops by one order
in case $c$, giving $(r,\ell)=(-2,-2m+2)$.
For \eqref{eq:def-second}, case $b$ uses the subleading left symbol
with $(r,\ell)=(-3,-2m+3)$, and case $c$ uses the subleading right
symbol with $(r,\ell)=(-2,-2m+2)$.
In each of Sections~\ref{sec:torsion}--\ref{sec:hodge}, 
$\Phi_1,\ldots,\Phi_5$ and $\widetilde\Phi_1,\ldots,\widetilde\Phi_5$ are the five local boundary densities in the first and second factorizations, respectively, in the order $a_1,a_2,a_3,b,c$.
In these local formulas the boundary densities and $d\operatorname{vol}_{\partial M}$
are understood at the chosen point $x_0\in\partial M$. The argument
$(x_0,\xi)$ in the integrand indicates evaluation after differentiation.
For example, the symbols at principal order in the first product are
$A i c(\xi)/|\xi|^2$ and $(|\xi|^2)^{1-m}$. The corresponding
symbols in the second are $A/|\xi|^2$ and
$i c(\xi)(|\xi|^2)^{1-m}$.
In taking a normal derivative one first keeps the coordinate covariable
fixed and uses $|\xi|^2=h(x_n)|\xi'|_{g^{\partial M}}^2+\xi_n^2$.
Only after differentiating do we set $x_n=0$ and
$|\xi'|_{g^{\partial M}}=1$. In the parallel normal frame this gives
\[
 \partial_{x_n}|\xi|^2=h'(0),\qquad
 \partial_{x_n}c(\xi')=\tfrac12 h'(0)c(\xi'),\qquad
 \partial_{x_n}c(dx_n)=0.
\]
The poles stay on their respective sides of the real axis for $h>0$.
A contour separating the poles can therefore be fixed locally in
$x_n$, which justifies commuting $\partial_{x_n}$ with the projection.

Two elementary integrals fix the normalizations used below. Closing in
the lower half-plane gives, for $p\geq1$ and $k\geq0$,
\begin{equation}
 \int_{\mathbb R}\frac{d\xi_n}{(\xi_n-i)^k(1+\xi_n^2)^p}
   =i^k\pi\,2^{2-k-2p}\binom{k+2p-2}{p-1}.
 \label{eq:normal-moment}
\end{equation}
Indeed the only enclosed pole is $-i$, of order $p$, and the contour
is clockwise. On the unit sphere in $T^*_{x_0}\partial M$,
\begin{equation}
 \int_{|\xi'|=1}\xi_a\,dS(\xi')=0,\qquad
 \int_{|\xi'|=1}\xi_a\xi_b\,dS(\xi')
   =\frac{\operatorname{Vol}(S^{2m-1})}{2m}\delta_{ab}.
 \label{eq:angular}
\end{equation}
All normal integrals use the real coordinate oriented from $-\infty$
to $+\infty$. The choice \eqref{eq:collar} fixes the signs of the
normal components and of $K$ throughout.

\section{The common interior residue}\label{sec:interior}

The odd-class parity argument \cite[Lemmas~1--2]{NO} applies to each of the three bundles. We give the critical symbol explicitly. The two choices of powers give the same pseudodifferential operator in the
interior. Indeed, let $\mathcal D$ denote any of the three operators and
let $c$ denote Clifford multiplication on its bundle $E$. Modulo
smoothing operators, the two products have the common value
\begin{align}
 &\bigl(c(X)c(Y)c(Z)\mathcal D^{-1}\bigr)\mathcal D^{-2m+2}
   =\bigl(c(X)c(Y)c(Z)\mathcal D^{-2}\bigr)\mathcal D^{-2m+3}
   \notag\\
 &\hspace{32mm}=c(X)c(Y)c(Z)\mathcal D^{-(2m-1)} .
 \label{int:product}
\end{align}

\begin{proposition}\label{int:vanishing}
For $D_T$ on $S(TM)$, $\widetilde D_F$ on $S(TM)\otimes F$, and
$D_{\mathrm{dR}}=d+\delta$ on $\Lambda^*T^*M\otimes\mathbb C$, the
interior residue density of \eqref{int:product} vanishes pointwise.
Consequently each of the six residues is the integral of its boundary density.
\end{proposition}

\begin{proof}
Set $L=\mathcal D^2$, so that $\mathcal D^{-(2m-1)}=\mathcal D L^{-m}$
modulo smoothing operators. Write $p_1+p_0$ for the symbol of
$\mathcal D$ and $q_{-2m-j}$ for the homogeneous component of degree
$-2m-j$ of the symbol of $L^{-m}$. The principal part $p_1=i c(\xi)$
is odd in $\xi$ and $p_0$ is independent of $\xi$. Since $L^m$ is a
differential operator of even order $2m$, its component of degree $2m-j$
has parity $(-1)^j$ under $\xi\mapsto-\xi$. Its parametrix starts with
$q_{-2m}=|\xi|^{-2m}$. At the next two orders the inverse symbol
recursion gives
\[
 q_{-2m-1}(x,-\xi)=-q_{-2m-1}(x,\xi),\qquad
 q_{-2m-2}(x,-\xi)=q_{-2m-2}(x,\xi).
\]
Indeed, in the recursion at order $-j$, a product with a symbol of
degree $2m-k$, a known inverse symbol of degree $-2m-l$, and
$|\alpha|$ covariable derivatives has $k+l+|\alpha|=j$. Its parity
is therefore $(-1)^j$.

The part of degree $-n=-2m-1$ in \eqref{int:product} is
\begin{equation}
 c(X)c(Y)c(Z)
 \left[
 p_1q_{-2m-2}+p_0q_{-2m-1}
 +\sum_{j=1}^{n}(\partial_{\xi_j}p_1)D_{x_j}q_{-2m-1}
 \right].
 \label{int:critical}
\end{equation}
The insertion multiplies from the left, hence no derivatives of $X,Y,Z$
enter this symbol. Each summand in the brackets is odd in $\xi$.
Taking the ordinary fiber trace preserves this parity. The interior
density at every $x\in M^\circ$ is thus
\[
 \left(\int_{|\xi|=1}\operatorname{tr}_E\!
 \left[\sigma_{-n}\!\left(c(X)c(Y)c(Z)
       \mathcal D^{-(2m-1)}\right)(x,\xi)\right]\,dS(\xi)\right)
 \,d\operatorname{vol}_M(x)=0 .
\]
In normal coordinates and a connection frame parallel at $x$, the three
choices of $p_0(x)$ are, respectively,
$\tfrac32c(T)-m c(V)$, $c(\omega)\otimes\Phi$, and $0$.
All three enter \eqref{int:critical} with the same covariable parity. The argument also applies to $d+\delta+\Psi$ for an arbitrary smooth bundle endomorphism $\Psi$, since it changes only $p_0$.
\end{proof}

\section{The Dirac operator with torsion}\label{sec:torsion}

The boundary symbol of \(D_T=D+\tfrac32c(T)-mc(V)\) contains both the
totally antisymmetric torsion and the vectorial torsion. We keep them separate
through the five cases. Denote the first five boundary densities by
\(\Phi_1,\ldots,\Phi_5\), in the order \(a_1,a_2,a_3,b,c\), and the
corresponding densities for the second factorization by
\(\widetilde\Phi_1,\ldots,\widetilde\Phi_5\).
In this section all traces are ordinary spinor traces.

\subsection{Boundary residue of \texorpdfstring{$\pi^+\!\bigl(c(X)c(Y)c(Z)D_T^{-1}\bigr)\circ\pi^+\!\bigl(D_T^{-2m+2}\bigr)$}{the torsion Dirac product with insertion c(X)c(Y)c(Z) and powers -1 and -2m+2}}\label{dt:subfirst}

\begin{theorem}\label{dt:first}
For \(m\geq2\) and arbitrary smooth vector fields \(X,Y,Z\),
\begin{equation}\label{dt:formula1}
\begin{aligned}
&\widetilde{\operatorname{Wres}}\!\left[\pi^+\!\left(c(X)c(Y)c(Z)D_T^{-1}\right)\circ\pi^+\!\left(D_T^{-2m+2}\right)\right]\\
&\quad=-\frac{\pi(2m-1)!\operatorname{Vol}(S^{2m-1})}{2^m(m+1)!(m-2)!}\int_{\partial M}\bigg\{
 \left.\partial_{x_n}\!\left[g(X,Y)Z_n-g(X,Z)Y_n+g(Y,Z)X_n\right]\right|_{\partial M}\\
&\qquad-\frac{m}{2}h'(0)[g(X,Y)Z_n-g(X,Z)Y_n+g(Y,Z)X_n]
 +\frac{6(m+1)}{2m-1}T(X,Y,Z)\bigg\}\,d\operatorname{vol}_{\partial M}.
\end{aligned}
\end{equation}
\end{theorem}

\begin{proof}
Fix \(x_0\in\partial M\). Each normal derivative is taken in the collar
before setting \(x_n=0\) and \(|\xi'|_{g^{\partial M}}=1\).
We evaluate the normal moments by \eqref{eq:normal-moment} and then use
\[
 \int_{S^{2m-1}}\xi_a\,dS=0,\qquad
 \int_{S^{2m-1}}\xi_a\xi_b\,dS=\frac{\operatorname{Vol}(S^{2m-1})}{2m}\delta_{ab}.
\]

\paragraph{Case \(a_1\).}
Take \(r=-1,\ \ell=-2m+2,\ j=k=0,\ |\alpha|=1\). The coefficient is \(-1\).
\begin{align*}
\Phi_1={}&-\int_{|\xi'|=1}\int_{\mathbb R}\sum_{a=1}^{2m}
 \operatorname{tr}_{S(TM)}\!\Bigl[
 \partial_{\xi_a}\pi^+_{\xi_n}\sigma_{-1}\!\left(c(X)c(Y)c(Z)D_T^{-1}\right)\\
&\qquad{}\cdot \partial_{x_a}\partial_{\xi_n}\sigma_{-2m+2}\!\left(D_T^{-2m+2}\right)\Bigr](x_0,\xi)
 \,d\xi_n\,dS(\xi')\,d\operatorname{vol}_{\partial M}.
\end{align*}
The principal symbols are \(iAc(\xi)/|\xi|^2\) and
\((|\xi|^2)^{1-m}\). The right factor
\(\partial_{x_a}\partial_{\xi_n}(|\xi|^2)^{1-m}\) vanishes for every
\(\xi_n\) since \(\partial_{x_a}|\xi|^2(x_0)=0\) for \(a\leq2m\).
Thus \(\Phi_1=0\).

\paragraph{Case \(a_2\).}
Take \(r=-1,\ \ell=-2m+2,\ j=1,\ k=|\alpha|=0\). The coefficient is \(-\frac12\).
\begin{align*}
\Phi_2={}&-\frac12\int_{|\xi'|=1}\int_{\mathbb R}
 \operatorname{tr}_{S(TM)}\!\Bigl[
 \partial_{x_n}\pi^+_{\xi_n}\sigma_{-1}\!\left(c(X)c(Y)c(Z)D_T^{-1}\right)\\
&\qquad{}\cdot \partial_{\xi_n}^{2}\sigma_{-2m+2}\!\left(D_T^{-2m+2}\right)\Bigr](x_0,\xi)
 \,d\xi_n\,dS(\xi')\,d\operatorname{vol}_{\partial M}.
\end{align*}
In the frame parallel along the normal, the upper-pole part of
\(ic(\xi)/|\xi|^2\) at \(x_0\) is
\([c(\xi')+i\,c(dx_n)]/[2(\xi_n-i)]\).
The pole moves with \(i\sqrt{h(x_n)}|\xi'|_{g^{\partial M}}\).
Differentiating the projection, including the insertion, gives
\begin{align*}
&\partial_{x_n}\pi^+_{\xi_n}\sigma_{-1}\!\left(c(X)c(Y)c(Z)D_T^{-1}\right)(x_0,\xi)\\
={}&\frac{\partial_{x_n}[c(X)c(Y)c(Z)]\,[c(\xi')+i\,c(dx_n)]}{2(\xi_n-i)}
+\frac{h'(0)c(X)c(Y)c(Z)[i\,c(\xi')-c(dx_n)]}{4(\xi_n-i)^2}.
\end{align*}
The right factor is
\(2(m-1)[(2m-1)\xi_n^2-1]/(1+\xi_n^2)^{m+1}\).
For \(j\le n\), Clifford anticommutation gives
\begin{equation}\label{dt:fourtrace}
 2^{-m}\operatorname{tr}\,[c(X)c(Y)c(Z)c(e_j)]
 =g(X,Y)Z_j-g(X,Z)Y_j+g(Y,Z)X_j,
\end{equation}
where \(c(e_n)=c(dx_n)\). Hence the trace at fixed \(\xi'\), before
multiplication by \(-1/2\), is
\begin{align*}
&\frac{2^{m+1}(m-1)[(2m-1)\xi_n^2-1]}{(1+\xi_n^2)^{m+1}}
 \bigg\{\frac{1}{2(\xi_n-i)}
 \bigg(\sum_{a=1}^{2m}\xi_a\,\partial_{x_n}[g(X,Y)Z_a-g(X,Z)Y_a+g(Y,Z)X_a]\\
&\qquad{}+i\,\partial_{x_n}[g(X,Y)Z_n-g(X,Z)Y_n+g(Y,Z)X_n]\bigg)\\
&\quad+\frac{h'(0)}{4(\xi_n-i)^2}
 \bigg(i\sum_{a=1}^{2m}\xi_a[g(X,Y)Z_a-g(X,Z)Y_a+g(Y,Z)X_a]\\
&\hspace{8em}{}-[g(X,Y)Z_n-g(X,Z)Y_n+g(Y,Z)X_n]\bigg)\bigg\}.
\end{align*}
The two normal moments are
\begin{align*}
\int_{\mathbb R}\frac{\partial_{\xi_n}^2(1+\xi_n^2)^{1-m}}{\xi_n-i}\,d\xi_n
 &=-\frac{4i\pi(2m-1)!}{2^{2m}(m+1)!(m-2)!},\\
\int_{\mathbb R}\frac{\partial_{\xi_n}^2(1+\xi_n^2)^{1-m}}{(\xi_n-i)^2}\,d\xi_n
 &=\frac{12m\pi(2m-1)!}{2^{2m}(m+2)!(m-2)!}.
\end{align*}
The integral of each term linear in \(\xi'\) vanishes over the tangent
sphere. Multiplication by \(-1/2\) gives
\begin{align*}
\Phi_2={}&-\frac{\pi(2m-1)!\operatorname{Vol}(S^{2m-1})}
 {2^m(m+1)!(m-2)!}\bigg\{
 \left.\partial_{x_n}[g(X,Y)Z_n-g(X,Z)Y_n+g(Y,Z)X_n]\right|_{\partial M}\\
&\quad 
-\frac{3m}{2(m+2)}h'(0)[g(X,Y)Z_n-g(X,Z)Y_n+g(Y,Z)X_n]
\bigg\}\,d\operatorname{vol}_{\partial M}.
\end{align*}

\paragraph{Case \(a_3\).}
Take \(r=-1,\ \ell=-2m+2,\ j=|\alpha|=0,\ k=1\). The coefficient is \(-\frac12\).
\begin{align*}
\Phi_3={}&-\frac12\int_{|\xi'|=1}\int_{\mathbb R}
 \operatorname{tr}_{S(TM)}\!\Bigl[
 \partial_{\xi_n}\pi^+_{\xi_n}\sigma_{-1}\!\left(c(X)c(Y)c(Z)D_T^{-1}\right)\\
&\qquad{}\cdot \partial_{\xi_n}\partial_{x_n}\sigma_{-2m+2}\!\left(D_T^{-2m+2}\right)\Bigr](x_0,\xi)
 \,d\xi_n\,dS(\xi')\,d\operatorname{vol}_{\partial M}.
\end{align*}
The differentiated projection is
\[
 -\frac{c(X)c(Y)c(Z)[c(\xi')+i\,c(dx_n)]}{2(\xi_n-i)^2}.
\]
The right factor is
\(\partial_{\xi_n}\partial_{x_n}(|\xi|^2)^{1-m}
=2m(m-1)h'(0)\xi_n/(1+\xi_n^2)^{m+1}\) at the boundary point.
The product has trace
\begin{align*}
&-\frac{2^m m(m-1)h'(0)\xi_n}{(\xi_n-i)^2(1+\xi_n^2)^{m+1}}
\bigg\{\sum_{a=1}^{2m}\xi_a[g(X,Y)Z_a-g(X,Z)Y_a+g(Y,Z)X_a]\\
&\quad  +i[g(X,Y)Z_n-g(X,Z)Y_n+g(Y,Z)X_n]\bigg\}.
\end{align*}
After multiplication by \(-1/2\), the required normal moment is
\[
\int_{\mathbb R}\frac{\xi_n\,d\xi_n}
 {(\xi_n-i)^2(1+\xi_n^2)^{m+1}}
 =\frac{i\pi(2m+2)!}{2^{2m+2}(m+1)!(m+2)!}.
\]
The constant term on the sphere gives
\begin{align*}
\Phi_3={}&-\frac{\pi(2m+1)!\operatorname{Vol}(S^{2m-1})}
 {2^{m+2}(m+2)!(m-2)!}\,
 h'(0)[g(X,Y)Z_n-g(X,Z)Y_n+g(Y,Z)X_n]\,d\operatorname{vol}_{\partial M}.
\end{align*}

\paragraph{Case \(b\).}
Take \(r=-1,\ \ell=-2m+1,\ j=k=|\alpha|=0\). The coefficient is \(-i\).
\begin{align*}
\Phi_4={}&-i\int_{|\xi'|=1}\int_{\mathbb R}
 \operatorname{tr}_{S(TM)}\!\Bigl[
 \pi^+_{\xi_n}\sigma_{-1}\!\left(c(X)c(Y)c(Z)D_T^{-1}\right)\\
&\qquad{}\cdot \partial_{\xi_n}\sigma_{-2m+1}\!\left(D_T^{-2m+2}\right)\Bigr](x_0,\xi)
 \,d\xi_n\,dS(\xi')\,d\operatorname{vol}_{\partial M}.
\end{align*}
Integration by parts in \(\xi_n\) moves the derivative to the
projected left principal symbol and changes the coefficient to \(+i\).
The derivative is the one used in case \(a_3\). Composing \(m-1\) factors
of \((D_T^2)^{-1}\) gives
\begin{align*}
\sigma_{-2m+1}(D_T^{-2m+2})
={}&-(m-1)(|\xi|^2)^{-m}\sigma_1(D_T^2)\\
&-\frac{im(m-1)}{2}(|\xi|^2)^{-m-1}
 \sum_{j=1}^{n}\partial_{\xi_j}(|\xi|^2)\,\partial_{x_j}(|\xi|^2).
\end{align*}
Substitution of the first-order symbol of \(D_T^2\) yields
\begin{align*}
\sigma_{-2m+1}(D_T^{-2m+2})(x_0,\xi)
={}&-i(m-1)\bigg\{
 \frac{h'(0)[m\xi_n+\tfrac12c(dx_n)c(\xi')]}{(1+\xi_n^2)^m}
 +\frac{m h'(0)\xi_n}{(1+\xi_n^2)^{m+1}}\\
&\quad 
 +\frac{c(\xi)(D_T-D)+(D_T-D)c(\xi)}{(1+\xi_n^2)^m}\bigg\}.
\end{align*}
For fixed \(\xi'\) the component trace of the term containing \(D_T-D\) is
\begin{align}
&2^{-m}\operatorname{tr}\!\left[c(X)c(Y)c(Z)c(e_j)
 \big(c(e_k)(D_T-D)+(D_T-D)c(e_k)\big)\right]\notag\\
&=2mV_k[g(X,Y)Z_j-g(X,Z)Y_j+g(Y,Z)X_j]\notag\\
&\quad-3[X_jT(e_k,Y,Z)-Y_jT(e_k,X,Z)+Z_jT(e_k,X,Y)]\notag\\
&\quad+3[g(X,Y)T(e_k,e_j,Z)-g(X,Z)T(e_k,e_j,Y)
             +g(Y,Z)T(e_k,e_j,X)].\label{dt:component}
\end{align}
The full trace is obtained by multiplying the right side
of \eqref{dt:component} by
\[
 \frac{2^m i(m-1)}{2(\xi_n-i)^2(1+\xi_n^2)^m}
 \bigg(\sum_{a=1}^{2m}\xi_a\delta_{aj}+i\delta_{nj}\bigg)\xi_k
\]
and summing over \(j,k\), then adding the metric term
\begin{align*}
&\frac{2^m i(m-1)h'(0)}{2(\xi_n-i)^2(1+\xi_n^2)^m}
 \left(m\xi_n-\frac i2+\frac{m\xi_n}{1+\xi_n^2}\right)\\
&\quad\times\left\{\sum_{a=1}^{2m}\xi_a[g(X,Y)Z_a-g(X,Z)Y_a+g(Y,Z)X_a]
 +i[g(X,Y)Z_n-g(X,Z)Y_n+g(Y,Z)X_n]\right\}.
\end{align*}
On multiplying by \(+i\), the basic normal integrals become
\begin{align*}
-\frac12\int_{\mathbb R}
\frac{d\xi_n}{(\xi_n-i)^2(1+\xi_n^2)^p}
 &=\frac{\pi(2p)!}{2^{2p+1}(p-1)!(p+1)!},\\
-\frac12\int_{\mathbb R}
\frac{\xi_n\,d\xi_n}{(\xi_n-i)^2(1+\xi_n^2)^p}
 &=-\frac{i\pi(2p)!}{2^{2p+1}p!(p+1)!}.
\end{align*}
Extracting \(\pi(2m)!/[2^{m+1}(m-2)!(m+1)!]\), the angularly
averaged diagonal weights in \eqref{dt:component} are \(1/(2m)\)
for \(j=k\leq2m\) and \(1/m\) for \(j=k=n\).
The off-diagonal terms disappear. This gives
\begin{align*}
\Phi_4={}&\frac{\pi(2m-1)!\operatorname{Vol}(S^{2m-1})}
 {2^{m+1}(m-2)!(m+1)!}
 \bigg\{\frac{m(5m+7)}{m+2}h'(0)
 [g(X,Y)Z_n-g(X,Z)Y_n+g(Y,Z)X_n]\\
&\quad-9T(X,Y,Z)
 -3[X_nT(\partial_{x_n},Y,Z)-Y_nT(\partial_{x_n},X,Z)
 +Z_nT(\partial_{x_n},X,Y)]\\
&\quad+2m[g(X,Y)g(Z,V)-g(X,Z)g(Y,V)+g(Y,Z)g(X,V)]\\
&\quad+2mV_n[g(X,Y)Z_n-g(X,Z)Y_n+g(Y,Z)X_n]\bigg\}
 \,d\operatorname{vol}_{\partial M}.
\end{align*}

\paragraph{Case \(c\).}
Take \(r=-2,\ \ell=-2m+2,\ j=k=|\alpha|=0\). The coefficient is \(-i\).
\begin{align*}
\Phi_5={}&-i\int_{|\xi'|=1}\int_{\mathbb R}
 \operatorname{tr}_{S(TM)}\!\Bigl[
 \pi^+_{\xi_n}\sigma_{-2}\!\left(c(X)c(Y)c(Z)D_T^{-1}\right)\\
&\qquad{}\cdot \partial_{\xi_n}\sigma_{-2m+2}\!\left(D_T^{-2m+2}\right)\Bigr](x_0,\xi)
 \,d\xi_n\,dS(\xi')\,d\operatorname{vol}_{\partial M}.
\end{align*}
The parametrix recursion, with
\(\sigma_0(D)(x_0)=-\frac m2h'(0)c(dx_n)\), gives
\begin{align*}
\sigma_{-2}(D_T^{-1})(x_0,\xi)
={}&-\frac{m h'(0)}2\frac{c(\xi)c(dx_n)c(\xi)}{(1+\xi_n^2)^2}
 +\frac{h'(0)}2\frac{c(\xi)c(dx_n)c(\xi')}{(1+\xi_n^2)^2}\\
&-h'(0)\frac{c(\xi)c(dx_n)c(\xi)}{(1+\xi_n^2)^3}
 +\frac{c(\xi)(D_T-D)c(\xi)}{(1+\xi_n^2)^2}.
\end{align*}
Taking the principal part at \(\xi_n=i\) gives
\begin{align*}
\pi^+_{\xi_n}\sigma_{-2}(D_T^{-1})
={}&-\frac{h'(0)[c(\xi')+i c(dx_n)]}{4(\xi_n-i)^3}
 +\frac{m h'(0)[c(dx_n)-i c(\xi')]}{4(\xi_n-i)^2}\\
&-\frac{[c(\xi')+i c(dx_n)](D_T-D)[c(\xi')+i c(dx_n)]}{4(\xi_n-i)^2}\\
&-\frac{i[c(\xi')(D_T-D)c(\xi')+c(dx_n)(D_T-D)c(dx_n)]}{4(\xi_n-i)}.
\end{align*}
The right factor is \(-2(m-1)\xi_n/(1+\xi_n^2)^m\).
For fixed \(\xi'\), multiplication by this factor and evaluation of the
ordinary trace with the insertion give
\begin{align*}
&\frac{2^{m-1}(m-1)\xi_n}{(1+\xi_n^2)^m}
\bigg\{
\bigg[\frac{h'(0)}{(\xi_n-i)^3}
 +\frac{im h'(0)}{(\xi_n-i)^2}
 +\frac{2m}{(\xi_n-i)^2}\Big(\sum_{a=1}^{2m}\xi_a V_a+iV_n\Big)\bigg]\\
&\quad\times\bigg[\sum_{b=1}^{2m}\xi_b
 \big(g(X,Y)Z_b-g(X,Z)Y_b+g(Y,Z)X_b\big)
 +i\big(g(X,Y)Z_n-g(X,Z)Y_n+g(Y,Z)X_n\big)\bigg]\\
&\quad-\frac{3}{(\xi_n-i)^2}
 \bigg[\Big(\sum_{a=1}^{2m}\xi_aX_a+iX_n\Big)
       \Big(\sum_{b=1}^{2m}\xi_bT(e_b,Y,Z)+iT(\partial_{x_n},Y,Z)\Big)\\
&\qquad\qquad-\Big(\sum_{a=1}^{2m}\xi_aY_a+iY_n\Big)
       \Big(\sum_{b=1}^{2m}\xi_bT(e_b,X,Z)+iT(\partial_{x_n},X,Z)\Big)\\
&\qquad\qquad+\Big(\sum_{a=1}^{2m}\xi_aZ_a+iZ_n\Big)
       \Big(\sum_{b=1}^{2m}\xi_bT(e_b,X,Y)+iT(\partial_{x_n},X,Y)\Big)\bigg]\\
&\quad+\frac{3i}{\xi_n-i}\bigg[T(X,Y,Z)
 -\sum_{a,b=1}^{2m}\xi_a\xi_b
 \big(X_aT(e_b,Y,Z)-Y_aT(e_b,X,Z)+Z_aT(e_b,X,Y)\big)\\
&\qquad\qquad-X_nT(\partial_{x_n},Y,Z)
 +Y_nT(\partial_{x_n},X,Z)-Z_nT(\partial_{x_n},X,Y)\bigg]\\
&\quad+\frac{2mi}{\xi_n-i}\bigg[
 \sum_{a,b=1}^{2m}\xi_a\xi_bV_a
 \big(g(X,Y)Z_b-g(X,Z)Y_b+g(Y,Z)X_b\big)\\
&\qquad\qquad+V_n\big(g(X,Y)Z_n-g(X,Z)Y_n+g(Y,Z)X_n\big)\\
&\qquad\qquad-g(X,Y)g(Z,V)+g(X,Z)g(Y,V)-g(Y,Z)g(X,V)\bigg]
\bigg\}.
\end{align*}
Here the terms involving \(T\) follow by commuting the two external
Clifford generators through \(c(T)\). After multiplying by \(-i\), the
normal moments for \(q=1,2,3\) are
\[
-i\,2^{m-1}(m-1)\int_{\mathbb R}
 \frac{\xi_n\,d\xi_n}{(1+\xi_n^2)^m(\xi_n-i)^q}
=\frac{\pi(2m-1)!}{2^{m+1}(m+1)!(m-2)!}
\begin{cases}
-2i(m+1)/(2m-1), & q=1,\\
2, & q=2,\\
3im/(m+2), & q=3.
\end{cases}
\]
Among the quadratic tangential moments only the diagonal ones survive. In particular,
\begin{align*}
&\sum_{a=1}^{2m}\big[X_aT(e_a,Y,Z)-Y_aT(e_a,X,Z)+Z_aT(e_a,X,Y)\big]\\
&\quad=3T(X,Y,Z)-X_nT(\partial_{x_n},Y,Z)
   +Y_nT(\partial_{x_n},X,Z)-Z_nT(\partial_{x_n},X,Y).
\end{align*}
Substitution yields
\begin{align*}
\Phi_5={}&\frac{\pi(2m-1)!\operatorname{Vol}(S^{2m-1})}
 {2^{m+1}(m+1)!(m-2)!}\bigg\{
 -\frac{m(2m+7)}{m+2}h'(0)
  \big[g(X,Y)Z_n-g(X,Z)Y_n+g(Y,Z)X_n\big]\\
&\quad+\frac{3(2m-7)}{2m-1}T(X,Y,Z)\\
&\quad+3\big[X_nT(\partial_{x_n},Y,Z)-Y_nT(\partial_{x_n},X,Z)
       +Z_nT(\partial_{x_n},X,Y)\big]\\
&\quad-2m\big[g(X,Y)g(Z,V)-g(X,Z)g(Y,V)+g(Y,Z)g(X,V)\big]\\
&\quad-2mV_n\big[g(X,Y)Z_n-g(X,Z)Y_n+g(Y,Z)X_n\big]\bigg\}
\,d\operatorname{vol}_{\partial M}.
\end{align*}

The \(V\)-terms in \(\Phi_4+\Phi_5\) cancel, as do the three separate terms
involving \(T(\partial_{x_n},\cdot,\cdot)\). The normal derivative comes
entirely from \(\Phi_2\). Since the interior density vanishes, integrating
the boundary density \(\phi_1=\sum_{j=1}^5\Phi_j\) gives \eqref{dt:formula1}.
\end{proof}

\subsection{Boundary residue of \texorpdfstring{$\pi^+\!\bigl(c(X)c(Y)c(Z)D_T^{-2}\bigr)\circ\pi^+\!\bigl(D_T^{-2m+3}\bigr)$}{the torsion Dirac product with insertion c(X)c(Y)c(Z) and powers -2 and -2m+3}}\label{dt:subsecond}

\begin{theorem}\label{dt:second}
Under the hypotheses of Theorem~\ref{dt:first}, the second factorization
has residue
\begin{equation}\label{dt:formula2}
\begin{aligned}
&\widetilde{\operatorname{Wres}}\!\left[\pi^+\!\left(c(X)c(Y)c(Z)D_T^{-2}\right)\circ\pi^+\!\left(D_T^{-2m+3}\right)\right]\\
&\quad=\frac{3\pi(2m-2)!\operatorname{Vol}(S^{2m-1})}{2^m(m+1)!(m-2)!}\int_{\partial M}\bigg\{
 \left.\partial_{x_n}\!\left[g(X,Y)Z_n-g(X,Z)Y_n+g(Y,Z)X_n\right]\right|_{\partial M}\\
&\qquad-\frac{m}{2}h'(0)[g(X,Y)Z_n-g(X,Z)Y_n+g(Y,Z)X_n]
 -2(m+1)T(X,Y,Z)\bigg\}\,d\operatorname{vol}_{\partial M}.
\end{aligned}
\end{equation}
\end{theorem}

\begin{proof}
The principal symbols of the two factors before truncation are
\(c(X)c(Y)c(Z)/|\xi|^2\) and \(ic(\xi)(|\xi|^2)^{1-m}\).

\paragraph{Case \(a_1\).}
Take \(r=-2,\ \ell=-2m+3,\ j=k=0,\ |\alpha|=1\). The coefficient is \(-1\).
\begin{align*}
\widetilde\Phi_1={}&-\int_{|\xi'|=1}\int_{\mathbb R}\sum_{a=1}^{2m}
 \operatorname{tr}_{S(TM)}\!\Bigl[
 \partial_{\xi_a}\pi^+_{\xi_n}\sigma_{-2}\!\left(c(X)c(Y)c(Z)D_T^{-2}\right)\\
&\qquad{}\cdot \partial_{x_a}\partial_{\xi_n}\sigma_{-2m+3}\!\left(D_T^{-2m+3}\right)\Bigr](x_0,\xi)
 \,d\xi_n\,dS(\xi')\,d\operatorname{vol}_{\partial M}.
\end{align*}
For each \(a\leq2m\), both \(\partial_{x_a}|\xi|^2\) and
\(\partial_{x_a}c(\xi)\) vanish at \(x_0\). Thus the mixed derivative
\[\partial_{\xi_n}\partial_{x_a}[ic(\xi)(|\xi|^2)^{1-m}]\] of the right
factor vanishes pointwise and \(\widetilde\Phi_1=0\).

\paragraph{Case \(a_2\).}
Take \(r=-2,\ \ell=-2m+3,\ j=1,\ k=|\alpha|=0\). The coefficient is \(-\frac12\).
\begin{align*}
\widetilde\Phi_2={}&-\frac12\int_{|\xi'|=1}\int_{\mathbb R}
 \operatorname{tr}_{S(TM)}\!\Bigl[
 \partial_{x_n}\pi^+_{\xi_n}\sigma_{-2}\!\left(c(X)c(Y)c(Z)D_T^{-2}\right)\\
&\qquad{}\cdot \partial_{\xi_n}^{2}\sigma_{-2m+3}\!\left(D_T^{-2m+3}\right)\Bigr](x_0,\xi)
 \,d\xi_n\,dS(\xi')\,d\operatorname{vol}_{\partial M}.
\end{align*}
The upper-pole part in the collar is
\[
 \frac{-i\,c(X)c(Y)c(Z)}
 {2\sqrt{h(x_n)}|\xi'|_{g^{\partial M}}
  [\xi_n-i\sqrt{h(x_n)}|\xi'|_{g^{\partial M}}]} .
\]
Differentiation of all three factors, followed by restriction to the
boundary and the tangential unit sphere, gives
\begin{align*}
&\partial_{x_n}\pi^+_{\xi_n}\sigma_{-2}\!\left(c(X)c(Y)c(Z)D_T^{-2}\right)\big|_{x_0,|\xi'|=1}\\
&\quad=-\frac{i}{2(\xi_n-i)}\,\partial_{x_n}[c(X)c(Y)c(Z)]
+h'(0)\left[\frac{i}{4(\xi_n-i)}+\frac{1}{4(\xi_n-i)^2}\right]c(X)c(Y)c(Z).
\end{align*}
The two Clifford components of the right factor are
\(i\,\partial_{\xi_n}^2[(1+\xi_n^2)^{1-m}]\) and
\(i\,\partial_{\xi_n}^2[\xi_n(1+\xi_n^2)^{1-m}]\).
By \eqref{dt:fourtrace} the fixed-\(\xi'\) trace is
\begin{align*}
2^m\sum_{j=1}^{n}
\partial_{\xi_n}^2\!\left[\frac{\xi_j}{(1+\xi_n^2)^{m-1}}\right]
\bigg\{&
\frac{\partial_{x_n}[g(X,Y)Z_j-g(X,Z)Y_j+g(Y,Z)X_j]}{2(\xi_n-i)}\\
&+h'(0)\left[-\frac1{4(\xi_n-i)}+\frac{i}{4(\xi_n-i)^2}\right]\\
&\quad {}\cdot[g(X,Y)Z_j-g(X,Z)Y_j+g(Y,Z)X_j]\bigg\}.
\end{align*}
For \(j=n\) the derivative also acts on the numerator \(\xi_n\).
Multiply by \(-1/2\) and integrate by parts twice. The corresponding
derivatives of the two projected coefficients are
\begin{align*}
\partial_{\xi_n}^2\frac1{2(\xi_n-i)}&=\frac1{(\xi_n-i)^3},\qquad
\partial_{\xi_n}^2\!\left[-\frac1{4(\xi_n-i)}+\frac{i}{4(\xi_n-i)^2}\right]
=-\frac1{2(\xi_n-i)^3}+\frac{3i}{2(\xi_n-i)^4}.
\end{align*}
Equation~\eqref{eq:normal-moment} gives the normal integral
\begin{align*}
\frac{3\pi(2m-2)!}{2^m(m-2)!(m+1)!}\bigg\{&
\partial_{x_n}[g(X,Y)Z_n-g(X,Z)Y_n+g(Y,Z)X_n]\\
&-\frac{5m}{2(m+2)}h'(0)[g(X,Y)Z_n-g(X,Z)Y_n+g(Y,Z)X_n]\\
&+\frac{i(2m-1)}{3}\sum_{a=1}^{2m}\xi_a
\Big(\partial_{x_n}[g(X,Y)Z_a-g(X,Z)Y_a+g(Y,Z)X_a]\\
&\quad {}-\frac{2m+1}{m+2}h'(0)[g(X,Y)Z_a-g(X,Z)Y_a+g(Y,Z)X_a]\Big)\bigg\}.
\end{align*}
The linear terms disappear on integration over \(S^{2m-1}\), giving
\begin{align*}
\widetilde\Phi_2
={}&\frac{3\pi(2m-2)!\operatorname{Vol}(S^{2m-1})}{2^m(m-2)!(m+1)!}
\bigg\{\left.\partial_{x_n}[g(X,Y)Z_n-g(X,Z)Y_n+g(Y,Z)X_n]
\right|_{\partial M}\\
&\quad {}-\frac{5m}{2(m+2)}h'(0)
[g(X,Y)Z_n-g(X,Z)Y_n+g(Y,Z)X_n]\bigg\}\,d\operatorname{vol}_{\partial M}.
\end{align*}

\paragraph{Case \(a_3\).}
Take \(r=-2,\ \ell=-2m+3,\ j=|\alpha|=0,\ k=1\). The coefficient is \(-\frac12\).
\begin{align*}
\widetilde\Phi_3={}&-\frac12\int_{|\xi'|=1}\int_{\mathbb R}
 \operatorname{tr}_{S(TM)}\!\Bigl[
 \partial_{\xi_n}\pi^+_{\xi_n}\sigma_{-2}\!\left(c(X)c(Y)c(Z)D_T^{-2}\right)\\
&\qquad{}\cdot \partial_{\xi_n}\partial_{x_n}\sigma_{-2m+3}\!\left(D_T^{-2m+3}\right)\Bigr](x_0,\xi)
 \,d\xi_n\,dS(\xi')\,d\operatorname{vol}_{\partial M}.
\end{align*}
The projected left symbol and its \(\xi_n\) derivative are
\(-iA/[2(\xi_n-i)]\) and \(iA/[2(\xi_n-i)^2]\).
After differentiation of the right principal symbol in the collar one finds
\begin{align*}
&\partial_{\xi_n}\partial_{x_n}\sigma_{-2m+3}(D_T^{-2m+3})\\
&\quad=\frac{i(m-1)h'(0)}{(1+\xi_n^2)^{m+1}}
\left[\xi_n(2m-1-\xi_n^2)c(\xi')+
\big((2m-1)\xi_n^2-1\big)c(dx_n)\right].
\end{align*}
The trace of the product is
\begin{align*}
&-\frac{2^{m-1}(m-1)h'(0)}{(\xi_n-i)^2(1+\xi_n^2)^{m+1}}\bigg\{
 \xi_n(2m-1-\xi_n^2)\sum_{a=1}^{2m}\xi_a
 [g(X,Y)Z_a-g(X,Z)Y_a+g(Y,Z)X_a]\\
&\quad {}+\big[(2m-1)\xi_n^2-1\big]
 [g(X,Y)Z_n-g(X,Z)Y_n+g(Y,Z)X_n]\bigg\}.
\end{align*}
After multiplication by \(-1/2\), the normal integral is
\begin{align*}
\frac{3\pi(2m)!\,h'(0)}{2^{m+2}(m-2)!(m+2)!}
\bigg\{&
g(X,Y)Z_n-g(X,Z)Y_n+g(Y,Z)X_n\\
&+\frac{2i(m^2-1)}{3m}\sum_{a=1}^{2m}\xi_a
[g(X,Y)Z_a-g(X,Z)Y_a+g(Y,Z)X_a]\bigg\}.
\end{align*}
The odd tangential term integrates to zero, so that
\begin{align*}
\widetilde\Phi_3={}&
\frac{3\pi(2m)!\operatorname{Vol}(S^{2m-1})}{2^{m+2}(m-2)!(m+2)!}\,
 h'(0)[g(X,Y)Z_n-g(X,Z)Y_n+g(Y,Z)X_n]\,d\operatorname{vol}_{\partial M}.
\end{align*}

\paragraph{Case \(b\).}
Take \(r=-3,\ \ell=-2m+3,\ j=k=|\alpha|=0\). The coefficient is \(-i\).
\begin{align*}
\widetilde\Phi_4={}&-i\int_{|\xi'|=1}\int_{\mathbb R}
 \operatorname{tr}_{S(TM)}\!\Bigl[
 \pi^+_{\xi_n}\sigma_{-3}\!\left(c(X)c(Y)c(Z)D_T^{-2}\right)\\
&\qquad{}\cdot \partial_{\xi_n}\sigma_{-2m+3}\!\left(D_T^{-2m+3}\right)\Bigr](x_0,\xi)
 \,d\xi_n\,dS(\xi')\,d\operatorname{vol}_{\partial M}.
\end{align*}
The left factor contains \(\sigma_{-3}(D_T^{-2})\). The parametrix recursion
gives
\begin{align*}
\sigma_{-3}(D_T^{-2})(x_0,\xi)
 &=-\frac{\sigma_1(D_T^2)(x_0,\xi)}{|\xi|^4}
   -\frac{2i\xi_n\,\partial_{x_n}(|\xi|^2)(x_0,\xi)}{|\xi|^6},\\
\sigma_1(D_T^2)(x_0,\xi)
 &=i\big[c(\xi)(D_T-D)+(D_T-D)c(\xi)\big]
   +im h'(0)\xi_n+\frac{i h'(0)}2c(dx_n)c(\xi').
\end{align*}
Its upper-pole part, including the insertion, is
\begin{align*}
&\pi^+_{\xi_n}\sigma_{-3}\!\left(c(X)c(Y)c(Z)D_T^{-2}\right)(x_0,\xi)\\
&=-i\,c(X)c(Y)c(Z)\bigg\{
 \left[-\frac1{4(\xi_n-i)^2}-\frac{i}{4(\xi_n-i)}\right]\!\\
&\quad \times\Big[c(\xi')(D_T-D)+(D_T-D)c(\xi')
       -\frac{h'(0)}2c(\xi')c(dx_n)\Big]\\
&\quad-\frac{i}{4(\xi_n-i)^2}
 \Big[c(dx_n)(D_T-D)+(D_T-D)c(dx_n)+m h'(0)\Big]
 -\frac{h'(0)}{4(\xi_n-i)^3}-\frac{i h'(0)}{8(\xi_n-i)^2}\bigg\}.
\end{align*}
The differentiated right principal symbol is
\[
\partial_{\xi_n}\sigma_{-2m+3}(D_T^{-2m+3})
=i\left[-\frac{2(m-1)\xi_n}{(1+\xi_n^2)^m}c(\xi')
 +\frac{1-(2m-3)\xi_n^2}{(1+\xi_n^2)^m}c(dx_n)\right].
\]
To write the trace at fixed \(\xi'\), the contribution containing \(D_T-D\)
has the component expansion
\begin{align*}
&2^m\sum_{j,k=1}^{n}
\pi^+_{\xi_n}\!\left[\frac{\xi_j}{(1+\xi_n^2)^2}\right]
\partial_{\xi_n}\!\left[\frac{\xi_k}{(1+\xi_n^2)^{m-1}}\right]
\Big\{2mV_j[g(X,Y)Z_k-g(X,Z)Y_k+g(Y,Z)X_k]\\
&\quad 
-3\big[X_kT(e_j,Y,Z)-Y_kT(e_j,X,Z)+Z_kT(e_j,X,Y)\\
&\quad 
+g(X,Y)T(e_j,e_k,Z)-g(X,Z)T(e_j,e_k,Y)+g(Y,Z)T(e_j,e_k,X)\big]\Big\}.
\end{align*}
The part proportional to \(h'(0)\) is
\begin{align*}
&2^m h'(0)[g(X,Y)Z_n-g(X,Z)Y_n+g(Y,Z)X_n]\,\bigg\{
 -\frac12\pi^+_{\xi_n}\!\left[\frac1{(1+\xi_n^2)^2}\right]
     \partial_{\xi_n}(1+\xi_n^2)^{1-m}\\
&\quad 
 +\pi^+_{\xi_n}\!\left[\frac{m\xi_n}{(1+\xi_n^2)^2}
                         +\frac{2\xi_n}{(1+\xi_n^2)^3}\right]
     \partial_{\xi_n}\!\left[\frac{\xi_n}{(1+\xi_n^2)^{m-1}}\right]\bigg\}\\
&{}+2^m h'(0)\sum_{a=1}^{2m}\xi_a[g(X,Y)Z_a-g(X,Z)Y_a+g(Y,Z)X_a]\,\bigg\{
 \frac12\pi^+_{\xi_n}\!\left[\frac1{(1+\xi_n^2)^2}\right]
     \partial_{\xi_n}\!\left[\frac{\xi_n}{(1+\xi_n^2)^{m-1}}\right]\\
&\quad 
 +\pi^+_{\xi_n}\!\left[\frac{m\xi_n}{(1+\xi_n^2)^2}
                         +\frac{2\xi_n}{(1+\xi_n^2)^3}\right]
     \partial_{\xi_n}(1+\xi_n^2)^{1-m}\bigg\}.
\end{align*}
Multiplying the trace by the boundary coefficient $-i$, we evaluate
the normal integrals by \eqref{eq:normal-moment}. The projected factor
$\pi^+_{\xi_n}[(1+\xi_n^2)^{-2}]$ gives
\begin{align*}
 -i\!\int_{\mathbb R}\!
 \pi^+_{\xi_n}\!\left[\frac{1}{(1+\xi_n^2)^2}\right]
 \partial_{\xi_n}\!\left[(1+\xi_n^2)^{1-m}\right]\,d\xi_n
 &=\frac{3m\pi(2m-2)!}{2^{2m}(m-2)!(m+1)!},\\
 -i\!\int_{\mathbb R}\!
 \pi^+_{\xi_n}\!\left[\frac{1}{(1+\xi_n^2)^2}\right]
 \partial_{\xi_n}\!\left[\xi_n(1+\xi_n^2)^{1-m}\right]\,d\xi_n
 &=-\frac{2i\pi(2m-1)!}{2^{2m}(m-1)!(m+1)!}.
\end{align*}
The two moments with $\pi^+_{\xi_n}[\xi_n(1+\xi_n^2)^{-2}]$ are
\begin{align*}
 -i\!\int_{\mathbb R}\!
 \pi^+_{\xi_n}\!\left[\frac{\xi_n}{(1+\xi_n^2)^2}\right]
 \partial_{\xi_n}\!\left[(1+\xi_n^2)^{1-m}\right]\,d\xi_n
 &=\frac{i\pi(2m-1)!}{2^{2m}(m-2)!(m+1)!},\\
 -i\!\int_{\mathbb R}\!
 \pi^+_{\xi_n}\!\left[\frac{\xi_n}{(1+\xi_n^2)^2}\right]
 \partial_{\xi_n}\!\left[\xi_n(1+\xi_n^2)^{1-m}\right]\,d\xi_n
 &=\frac{3\pi(2m-2)!}{2^{2m}(m-2)!(m+1)!}.
\end{align*}
For the projected term with denominator $(1+\xi_n^2)^3$ we obtain
\begin{align*}
 -i\!\int_{\mathbb R}\!
 \pi^+_{\xi_n}\!\left[\frac{\xi_n}{(1+\xi_n^2)^3}\right]
 \partial_{\xi_n}\!\left[(1+\xi_n^2)^{1-m}\right]\,d\xi_n
 &=\frac{i\pi(2m-1)!(2m+1)}{2^{2m+1}(m-2)!(m+2)!},\\
 -i\!\int_{\mathbb R}\!
 \pi^+_{\xi_n}\!\left[\frac{\xi_n}{(1+\xi_n^2)^3}\right]
 \partial_{\xi_n}\!\left[\xi_n(1+\xi_n^2)^{1-m}\right]\,d\xi_n
 &=\frac{15m\pi(2m-2)!}{2^{2m+2}(m-2)!(m+2)!}.
\end{align*}
For example, integrating by parts in the moment with
$\pi^+_{\xi_n}[\xi_n(1+\xi_n^2)^{-2}]$ and
$\partial_{\xi_n}[\xi_n(1+\xi_n^2)^{1-m}]$ gives
\[
-\frac12\int_{\mathbb R}
 \frac{\xi_n\,d\xi_n}{(\xi_n-i)^3(1+\xi_n^2)^{m-1}}
=\frac{3\pi(2m-2)!}{2^{2m}(m-2)!(m+1)!}.
\]
After factoring out $3\pi(2m-2)!/[2^{2m}(m-2)!(m+1)!]$ and taking the angular average, the tangential diagonal terms have weight
\(1/2\), and the normal diagonal term has weight \(1\).
Consequently
\begin{align*}
\widetilde\Phi_4
={}&\frac{3\pi(2m-2)!\operatorname{Vol}(S^{2m-1})}
 {2^m(m-2)!(m+1)!}
 \bigg\{\left[\frac{m(m+7)}{2(m+2)}h'(0)+mV_n\right]
 [g(X,Y)Z_n-g(X,Z)Y_n+g(Y,Z)X_n]\\
&\quad-\frac92 T(X,Y,Z)
 -\frac32\big[X_nT(\partial_{x_n},Y,Z)-Y_nT(\partial_{x_n},X,Z)
                          +Z_nT(\partial_{x_n},X,Y)\big]\\
&\quad+m\big[g(X,Y)g(Z,V)-g(X,Z)g(Y,V)+g(Y,Z)g(X,V)\big]
 \bigg\}\,d\operatorname{vol}_{\partial M}.
\end{align*}

\paragraph{Case \(c\).}
Take \(r=-2,\ \ell=-2m+2,\ j=k=|\alpha|=0\). The coefficient is \(-i\).
\begin{align*}
\widetilde\Phi_5={}&-i\int_{|\xi'|=1}\int_{\mathbb R}
 \operatorname{tr}_{S(TM)}\!\Bigl[
 \pi^+_{\xi_n}\sigma_{-2}\!\left(c(X)c(Y)c(Z)D_T^{-2}\right)\\
&\qquad{}\cdot \partial_{\xi_n}\sigma_{-2m+2}\!\left(D_T^{-2m+3}\right)\Bigr](x_0,\xi)
 \,d\xi_n\,dS(\xi')\,d\operatorname{vol}_{\partial M}.
\end{align*}
Move the \(\xi_n\) derivative to the projected left symbol by integration
by parts. The resulting factor is \(+i\) times
\(iA/[2(\xi_n-i)^2]\). The right subleading symbol is obtained from
\((D_T^{-2})^{m-1}D_T\). In the boundary chart its symbol of the required degree is
\begin{align*}
\sigma_{-2m+2}(D_T^{-2m+3})
={}&(|\xi|^2)^{1-m}\sigma_0(D_T)
  +\sigma_{-2m+1}\!\left((D_T^{-2})^{m-1}\right)i\,c(\xi)\\
 &-i\,\partial_{\xi_n}\!\left[(|\xi|^2)^{1-m}\right]
            \partial_{x_n}[i\,c(\xi)],\\
\sigma_{-2m+1}\!\left((D_T^{-2})^{m-1}\right)
={}&(m-1)(|\xi|^2)^{2-m}\sigma_{-3}(D_T^{-2})\\
&-i\binom{m-1}{2}(|\xi|^2)^{3-m}
      \partial_{\xi_n}(|\xi|^{-2})\,\partial_{x_n}(|\xi|^{-2}).
\end{align*}
Using the first normal derivatives of \(|\xi|^2\) and \(c(\xi)\), we find
\begin{align*}
&\sigma_{-2m+2}(D_T^{-2m+3})(x_0,\xi)\\
={}&h'(0)\left[
 \frac{(m-1)(m-\tfrac12)\xi_n}{(1+\xi_n^2)^m}
 +\frac{m(m-1)\xi_n}{(1+\xi_n^2)^{m+1}}
 \right]c(\xi')\\
&+h'(0)\left[
 -\frac{m}{2(1+\xi_n^2)^{m-1}}
 +\frac{(m-1)(m\xi_n^2-\tfrac12)}{(1+\xi_n^2)^m}
 +\frac{m(m-1)\xi_n^2}{(1+\xi_n^2)^{m+1}}
 \right]c(dx_n)\\
&+\frac{D_T-D}{(1+\xi_n^2)^{m-1}}
 +\frac{m-1}{(1+\xi_n^2)^m}
 \big[c(\xi)(D_T-D)+(D_T-D)c(\xi)\big]c(\xi).
\end{align*}
The trace at fixed \(\xi'\), including the coefficient after integration
by parts, is
\begin{align*}
&-\frac{2^{m-1}}{(\xi_n-i)^2}\bigg\{\!
 h'(0)\left[\frac{(m-1)(m-\tfrac12)\xi_n}{(1+\xi_n^2)^m}
      +\frac{m(m-1)\xi_n}{(1+\xi_n^2)^{m+1}}\right]\\[-1mm]
&\qquad\times\sum_{a=1}^{2m}\xi_a
 [g(X,Y)Z_a-g(X,Z)Y_a+g(Y,Z)X_a]\\
&\quad+h'(0)\left[-\frac{m}{2(1+\xi_n^2)^{m-1}}
 +\frac{(m-1)(m\xi_n^2-\tfrac12)}{(1+\xi_n^2)^m}
 +\frac{m(m-1)\xi_n^2}{(1+\xi_n^2)^{m+1}}\right]\\[-1mm]
&\qquad\times[g(X,Y)Z_n-g(X,Z)Y_n+g(Y,Z)X_n]\\
&\quad+\frac{1}{(1+\xi_n^2)^{m-1}}\Big[\tfrac32T(X,Y,Z)
 -m\big(g(X,Y)g(Z,V)-g(X,Z)g(Y,V)+g(Y,Z)g(X,V)\big)\Big]\\
&\quad+\frac{m-1}{(1+\xi_n^2)^m}\sum_{j,k=1}^n\xi_j\xi_k\Big[
 2mV_j\big(g(X,Y)Z_k-g(X,Z)Y_k+g(Y,Z)X_k\big)\\
&\quad -3\big(X_kT(e_j,Y,Z)-Y_kT(e_j,X,Z)
 +Z_kT(e_j,X,Y)\big)\Big]\bigg\}.
\end{align*}
The terms containing \(T(e_j,e_k,\cdot)\) vanish in the symmetric sum in
\(j,k\). Besides the first two moments used for the first factorization,
the normal integral requires
\[
 -\frac12\int_{\mathbb R}
 \frac{\xi_n^2\,d\xi_n}{(\xi_n-i)^2(1+\xi_n^2)^p}
 =\frac{\pi(2p-2)!(p-2)}{2^{2p}(p-1)!(p+1)!},\qquad p\geq1.
\]
Extract \(\pi(2m-2)!/[2^m(m-2)!(m+1)!]\).
For \(j,k\leq2m\) the last component sum in the trace then has weight
\(m(2m-1)\xi_j\xi_k\). For \(j=k=n\) the weight is \(m-2\).
The cross terms are odd in \(\xi'\).
The tangential and normal diagonal coefficients after averaging are
therefore \((2m-1)/2\) and \(m-2\). Substitution gives
\begin{align*}
\widetilde\Phi_5
={}&\frac{3\pi(2m-2)!\operatorname{Vol}(S^{2m-1})}
 {2^m(m-2)!(m+1)!}\bigg\{
 -\left[\frac{m(4m+3)}{2(m+2)}h'(0)+mV_n\right]\\
&\quad\times [g(X,Y)Z_n-g(X,Z)Y_n+g(Y,Z)X_n]
 +\frac{5-4m}{2}T(X,Y,Z)\\
&\quad+\frac32\big[X_nT(\partial_{x_n},Y,Z)
 -Y_nT(\partial_{x_n},X,Z)+Z_nT(\partial_{x_n},X,Y)\big]\\
&\quad-m\big[g(X,Y)g(Z,V)-g(X,Z)g(Y,V)+g(Y,Z)g(X,V)\big]
\bigg\}\,d\operatorname{vol}_{\partial M}.
\end{align*}

In \(\widetilde\Phi_4+\widetilde\Phi_5\) the \(V\)-terms and the
separate terms involving \(T(\partial_{x_n},\cdot,\cdot)\) cancel.
With the prefactor in \eqref{dt:formula2}, the metric and torsion
coefficients are \(-m/2\) and \(-2(m+1)\), respectively. Integrating the
boundary density \(\phi_2=\sum_{j=1}^5\widetilde\Phi_j\) proves \eqref{dt:formula2}.
\end{proof}

\subsection{Separating skew torsion from mean curvature}\label{dt:interpretation}

Each quantity in square brackets in \eqref{dt:formula1} and
\eqref{dt:formula2} is a contraction with the collar normal. Its normal
derivative differentiates the entire scalar before restriction and is
independent of the tangential frame. Thus both integrands define global
boundary scalars for the given fields on the collar. Multiplication by
\(d\operatorname{vol}_{\partial M}\) gives their densities.
The common coefficient of \(T(X,Y,Z)\) allows the two residues to separate
extrinsic geometry from torsion.

\begin{corollary}\label{dt:separation}
For arbitrary \(X,Y,Z\) the weighted sum satisfies
\begin{align}
 &\frac{3}{2m+2}\widetilde{\operatorname{Wres}}\!\left[\pi^+\!\left(c(X)c(Y)c(Z)D_T^{-1}\right)\circ\pi^+\!\left(D_T^{-2m+2}\right)\right]\notag\\
 &\quad{}+\frac{2m-1}{2m+2}\widetilde{\operatorname{Wres}}\!\left[\pi^+\!\left(c(X)c(Y)c(Z)D_T^{-2}\right)\circ\pi^+\!\left(D_T^{-2m+3}\right)\right]\notag\\
 &\qquad=-\frac{6\pi(2m-2)!\operatorname{Vol}(S^{2m-1})}
                  {2^m m!(m-2)!}
 \int_{\partial M}T(X,Y,Z)\,d\operatorname{vol}_{\partial M}.\label{dt:pureT}
\end{align}
If \(g(X,Y)Z_n-g(X,Z)Y_n+g(Y,Z)X_n=C_0\) is constant throughout a
collar, then
\begin{align}
 &\widetilde{\operatorname{Wres}}\!\left[\pi^+\!\left(c(X)c(Y)c(Z)D_T^{-2}\right)\circ\pi^+\!\left(D_T^{-2m+3}\right)\right]\notag\\
 &\quad{}-\widetilde{\operatorname{Wres}}\!\left[\pi^+\!\left(c(X)c(Y)c(Z)D_T^{-1}\right)\circ\pi^+\!\left(D_T^{-2m+2}\right)\right]\notag\\
 &\qquad=\frac{\pi(2m-2)!\operatorname{Vol}(S^{2m-1})}{2^m m!(m-2)!}
 C_0\int_{\partial M}K\,d\operatorname{vol}_{\partial M}.\label{dt:pureK}
\end{align}
\end{corollary}
\begin{proof}
Use \(K=-mh'(0)\). The weights in \eqref{dt:pureT} cancel both geometric
terms and add to one. Under the collar hypothesis the normal derivative
in each formula vanishes. Subtracting then gives \eqref{dt:pureK}.
\end{proof}

Taking \(X=Y=Z=\partial_{x_n}\) throughout a smaller collar yields
\(C_0=1\) and \(T(X,Y,Z)=0\). Such fields extend smoothly to \(M\) after
multiplication by a cutoff equal to one near the boundary.
Formula~\eqref{dt:pureK} therefore realizes the integrated trace of the
second fundamental form for every torsion in this family.
This is the geometric integral entering the Riemannian gravitational
boundary action \cite{GH}, with our inward normal and the numerical normalization in
\eqref{dt:pureK}. The formulas depend linearly on the skew torsion
sampled by the three fields. In dimension five the coefficients of the
geometric bracket in the two residues are \(-\pi^3/2\) and \(+\pi^3/2\),
respectively, and the torsion coefficient is \(-3\pi^3\) in both.

\begin{corollary}\label{dt:tangency}
Suppose \(X_n=Y_n=Z_n=0\) on \(\partial M\). The \(h'(0)\) terms in
\eqref{dt:formula1} and \eqref{dt:formula2} vanish and the first-order
term reduces to
\[
 \left.[g(X,Y)\partial_{x_n}Z_n-g(X,Z)\partial_{x_n}Y_n
             +g(Y,Z)\partial_{x_n}X_n]\right|_{\partial M}.
\]
If \(X,Y,Z\) are tangent to every slice of a collar, the two residues
coincide, and each is the right-hand side of \eqref{dt:pureT}.
\end{corollary}
\begin{proof}
Apply the product rule to the complete bracket in the two theorems.
Tangency along the collar forces the normal components to vanish identically.
\end{proof}

There is also an intrinsic interpretation on the closed even-dimensional
spin manifold \(\partial M\), in terms of the spectral torsion functional of \cite{DSZ}.
Give the boundary its induced spin structure. Write
\(D^{\partial M}\) for its spin Dirac operator and \(c_{\partial M}\) for
its Clifford action. Restrict the three-form \(T\) to tangent arguments
and set
\[
 B=D^{\partial M}+\tfrac32c_{\partial M}(T|_{T\partial M}).
\]
For tangent boundary fields \(u,v,w\), let \(X,Y,Z\) be any smooth extensions.
With the same residue normalization, the weighted sum on the left of
\eqref{dt:pureT} equals
\begin{equation}\label{dt:intrinsic}
 \frac{\pi(2m-2)!}{2^{2m-1}m!(m-2)!}
 \operatorname{Wres}_{\partial M}\!\left[
 c_{\partial M}(u)c_{\partial M}(v)c_{\partial M}(w)B^{-(2m-1)}\right].
\end{equation}
To see this, use geodesic coordinates and a synchronous spin frame on
\(\partial M\), and put \(P=\tfrac32c_{\partial M}(T|_{T\partial M})\).
At the center of the frame the degree \(-2m\) symbol of \(B(B^2)^{-m}\)
on the unit cotangent sphere is
\[
 (1-m)P+m c_{\partial M}(\eta)P c_{\partial M}(\eta).
\]
Its angular average is \(-2P\), since
\[\sum_{a=1}^{2m}c_{\partial M}(e_a)P c_{\partial M}(e_a)=(2m-6)P\].
The trace of
\[c_{\partial M}(u)c_{\partial M}(v)c_{\partial M}(w)c_{\partial M}(T|_{T\partial M})\]
is \(2^mT(u,v,w)\). Thus the intrinsic residue equals
\[-3\cdot2^m\operatorname{Vol}(S^{2m-1})
 \int_{\partial M}T(u,v,w)\,d\operatorname{vol}_{\partial M}.\]
Substitution in \eqref{dt:pureT} proves \eqref{dt:intrinsic} for all
extensions. This identifies the boundary torsion sampled by the two
factorizations with the torsion of a Dirac operator on the boundary itself.

\paragraph{Local reconstruction and the normal component of torsion.}
Let $\mathcal S_T(X,Y,Z)$ denote the weighted functional in \eqref{dt:pureT}, and let $b_m$ be its displayed nonzero coefficient. For a smooth boundary function $f$, extended to the collar in any manner, that identity gives
\[
 \mathcal S_T(fX,Y,Z)=b_m\int_{\partial M}f\,T(X,Y,Z)\,d\operatorname{vol}_{\partial M}.
\]
Thus arbitrary localized tests determine the smooth scalar $T(X,Y,Z)$ pointwise. Choosing local frames and using cutoffs recovers the restriction of the ambient tensor to $TM|_{\partial M}$. This is stronger than recovering only the pullback of $T$ to the boundary. Indeed, with $\nu=\partial_{x_n}$, write at the boundary
\[
 T=\tau+\nu^\flat\wedge\beta,\qquad
 \tau=T|_{\Lambda^3T\partial M},\qquad
 \beta=(\iota_\nu T)|_{\Lambda^2T\partial M}.
\]
Tangential triples determine $\tau$, while $(\nu,u,v)$ determines $\beta(u,v)$. An intrinsic boundary Dirac operator detects the former component in \eqref{dt:intrinsic}; normal insertions supply additional ambient information.

\paragraph{An extension-independent torsion response.}
The preceding localization works because the weighted combination cancels the entire first normal derivative, including derivatives of a multiplying test function. Individual residues need not descend to a functional of boundary fields alone. The cancellation therefore has a geometric meaning beyond simplifying coefficients: it removes the choice of continuation of the probes off the boundary. Conversely, the difference of the residues isolates their common first-order geometric expression. With normally parallel probes its derivative part vanishes, leaving the mean-curvature contribution. These are two complementary measurements, one of an alternating coupling and the other of the normal change of geometry.

\paragraph{Interpretation as a probe of a torsion coupling.}
In a fermionic differential expression the term $\tfrac32c(T)$ is a local coupling to the three-form background. Formula~\eqref{dt:pureT} supplies a linear, localized boundary probe of that coupling. This is compatible with the role of skew torsion in spinorial geometry \cite{Agricola,PS}; it does not require identifying $T$ with a particular physical field. If an application identifies it with a three-form field strength, closedness, flux quantization, and field equations are additional hypotheses. No such hypotheses follow from the residue identity. Likewise, cancellation of vectorial torsion in the final answer means that these observables do not distinguish that component; it does not mean that the underlying Dirac operator or other spectral quantities are independent of it.

\section{The Connes-type operator on twisted spinors}\label{sec:connes}
Wang and Wang define a Connes-type operator on an even-dimensional spin
manifold in \cite[Definition~2.2, Eq.~(2.8)]{WW}. We use their differential
expression on the odd-dimensional bundle $S(TM)\otimes F$,
\[
 \widetilde D_F=D_F+c(\omega)\otimes\Phi,\qquad
 D_F=\sum_{j=1}^{2m+1}c(e_j)\nabla^{S(TM)\otimes F}_{e_j}.
\]
Here $F$ is a Hermitian vector bundle with unitary connection,
$\Phi=\Phi^*$ is a smooth endomorphism of $F$, and $\omega$ is a smooth
complex differential form, possibly of mixed degree. We let $A=c(X)c(Y)c(Z)$
act on the spinor factor. Two distinct fiber traces occur. A term without
the perturbation acquires the multiplicity $\operatorname{rk}(F)$.
A term linear in the perturbation acquires $\operatorname{tr}_F(\Phi)$.
The order constraint at the boundary allows precisely one
subleading symbol or one first-order loss from differentiation. Thus terms
quadratic in $\Phi$ and derivatives of $\Phi$ cannot enter.

\subsection{Boundary residue of \texorpdfstring{$\pi^+\!\bigl(c(X)c(Y)c(Z)\widetilde D_F^{-1}\bigr)\circ\pi^+\!\bigl(\widetilde D_F^{-2m+2}\bigr)$}{the Connes-type product with insertion c(X)c(Y)c(Z) and powers -1 and -2m+2}}
In an oriented orthonormal frame fix the irreducible spin representation by
$c(e_1)\cdots c(e_n)=\kappa_n=i^{-(m+1)}$. Let $*$ denote the
complex-linear Hodge star of $M$. For an ordered unit simple $q$-form
$\beta$, the identities $c(\beta)c(*\beta)=\kappa_n$ and
$c(\beta)^2=(-1)^{q(q+1)/2}$ imply
$c(*\beta)=\kappa_n(-1)^{q(q+1)/2}c(\beta)$. By linearity the last
identity holds for every homogeneous $\beta$. Consequently $c(\omega)$
can be represented by odd-degree forms alone. The components of this
representative that can occur in a trace with at most five other
Clifford factors are
\begin{align}
 \theta_1&=\omega_{[1]}-\kappa_n*\omega_{[n-1]},&
 \theta_3&=\omega_{[3]}+\kappa_n*\omega_{[n-3]},&
 \theta_5&=\omega_{[5]}-\kappa_n*\omega_{[n-5]}.\label{co:dual}
\end{align}
Here $\omega_{[q]}$ is the degree-$q$ component. This reduction matters
already in dimension five, where a two-form contributes through its Hodge
dual.

\begin{theorem}\label{co:one}
For arbitrary smooth $X,Y,Z$ and $m\ge2$, the first residue of
$\widetilde D_F$ equals
\begin{equation}\label{co:first}
\begin{aligned}
 &\widetilde{\operatorname{Wres}}\!\left[\pi^+\!\left(c(X)c(Y)c(Z)\widetilde D_F^{-1}\right)\circ\pi^+\!\left(\widetilde D_F^{-2m+2}\right)\right]\\
 &\quad=-\frac{\pi(2m-1)!\operatorname{Vol}(S^{2m-1})}
              {2^m(m+1)!(m-2)!}\int_{\partial M}\!\bigg\{\\
 &\qquad\operatorname{rk}(F)
 \left.\partial_{x_n}\big[g(X,Y)Z_n-g(X,Z)Y_n+g(Y,Z)X_n\big]
 \right|_{x_n=0}\\
 &\quad-\frac m2\operatorname{rk}(F)h'(0)
          [g(X,Y)Z_n-g(X,Z)Y_n+g(Y,Z)X_n]\\
 &\quad+\frac{4(m+1)}{2m-1}\operatorname{tr}_F(\Phi)
       (\omega_{[3]}+\kappa_n*\omega_{[n-3]})(X,Y,Z)\bigg\}
             \,d\operatorname{vol}_{\partial M}.
\end{aligned}
\end{equation}
\end{theorem}

\begin{proof}
Choose a unitary frame of $F$ with vanishing connection matrix at the
boundary point $x_0$. The principal symbol is $ic(\xi)$ and the order-zero
symbol is $-mh'(0)c(dx_n)/2+c(\omega)\otimes\Phi$. In particular, the
metric parts of all five cases are the untwisted spin contributions of
Section~\ref{sec:torsion} multiplied by $\operatorname{rk}(F)$.
\paragraph{Case \(a_1\).}
Take \(r=-1,\ \ell=-2m+2,\ j=k=0,\ |\alpha|=1\). The coefficient is \(-1\).
\begin{align*}
\Phi_1={}&-\int_{|\xi'|=1}\int_{\mathbb R}\sum_{a=1}^{2m}
 \operatorname{tr}_{S(TM)\otimes F}\!\Bigl[
 \partial_{\xi_a}\pi^+_{\xi_n}\sigma_{-1}\!\left(c(X)c(Y)c(Z)\widetilde D_F^{-1}\right)\\
&\qquad{}\cdot \partial_{x_a}\partial_{\xi_n}\sigma_{-2m+2}\!\left(\widetilde D_F^{-2m+2}\right)\Bigr](x_0,\xi)
 \,d\xi_n\,dS(\xi')\,d\operatorname{vol}_{\partial M}.
\end{align*}
In tangential geodesic coordinates $\partial_{x_a}|\xi|^2(x_0)=0$.
Thus $\partial_{x_a}\partial_{\xi_n}(|\xi|^2)^{1-m}=0$ at $x_0$ for
every $a\leq 2m$ and every $\xi_n$. The integrand vanishes pointwise, giving $\Phi_1=0$.

\paragraph{Case \(a_2\).}
Take \(r=-1,\ \ell=-2m+2,\ j=1,\ k=|\alpha|=0\). The coefficient is \(-\frac12\).
\begin{align*}
\Phi_2={}&-\frac12\int_{|\xi'|=1}\int_{\mathbb R}
 \operatorname{tr}_{S(TM)\otimes F}\!\Bigl[
 \partial_{x_n}\pi^+_{\xi_n}\sigma_{-1}\!\left(c(X)c(Y)c(Z)\widetilde D_F^{-1}\right)\\
&\qquad{}\cdot \partial_{\xi_n}^{2}\sigma_{-2m+2}\!\left(\widetilde D_F^{-2m+2}\right)\Bigr](x_0,\xi)
 \,d\xi_n\,dS(\xi')\,d\operatorname{vol}_{\partial M}.
\end{align*}
The moving upper pole is
$i\sqrt{h(x_n)}|\xi'|_{g^{\partial M}}$. Differentiating the projected
left symbol along the collar, including the insertion, and then taking
$x_n=0$ and $|\xi'|_{g^{\partial M}}=1$ gives
\[
 \frac{\partial_{x_n}[c(X)c(Y)c(Z)]\,[c(\xi')+i c(dx_n)]}
      {2(\xi_n-i)}
 +\frac{h'(0)A[i c(\xi')-c(dx_n)]}{4(\xi_n-i)^2}.
\]
The right factor is $2(m-1)[(2m-1)\xi_n^2-1]/(1+\xi_n^2)^{m+1}$.
The terms with $c(\xi')$ have zero angular average.
The trace of $A c(dx_n)$ is
$2^m\operatorname{rk}(F)[g(X,Y)Z_n-g(X,Z)Y_n+g(Y,Z)X_n]$
on $S(TM)\otimes F$. Its normal derivative gives the derivative of the
entire bracket. The two normal moments in Section~\ref{dt:subfirst},
case $a_2$, give
\begin{align*}
 &\Phi_2=-\frac{\operatorname{rk}(F)\,\pi(2m-1)!\operatorname{Vol}(S^{2m-1})}{2^{m}(m+1)!(m-2)!}\bigg\{\\
 &\hspace{4em}\left.\partial_{x_n}\!\left[g(X,Y)Z_n-g(X,Z)Y_n+g(Y,Z)X_n\right]\right|_{x_n=0}\\
 &\hspace{4em}-\frac{3m}{2(m+2)}h'(0)[g(X,Y)Z_n-g(X,Z)Y_n+g(Y,Z)X_n]\bigg\}\,d\operatorname{vol}_{\partial M}.
\end{align*}

\paragraph{Case \(a_3\).}
Take \(r=-1,\ \ell=-2m+2,\ j=|\alpha|=0,\ k=1\). The coefficient is \(-\frac12\).
\begin{align*}
\Phi_3={}&-\frac12\int_{|\xi'|=1}\int_{\mathbb R}
 \operatorname{tr}_{S(TM)\otimes F}\!\Bigl[
 \partial_{\xi_n}\pi^+_{\xi_n}\sigma_{-1}\!\left(c(X)c(Y)c(Z)\widetilde D_F^{-1}\right)\\
&\qquad{}\cdot \partial_{\xi_n}\partial_{x_n}\sigma_{-2m+2}\!\left(\widetilde D_F^{-2m+2}\right)\Bigr](x_0,\xi)
 \,d\xi_n\,dS(\xi')\,d\operatorname{vol}_{\partial M}.
\end{align*}
At $x_0$, the two factors after differentiation are
\[
 -\frac{A[c(\xi')+i c(dx_n)]}{2(\xi_n-i)^2},
 \qquad
 \frac{2m(m-1)h'(0)\xi_n}{(1+\xi_n^2)^{m+1}} .
\]
The first moment on the tangential sphere removes $c(\xi')$.
In the normal component, integration of
$\xi_n/[(\xi_n-i)^2(1+\xi_n^2)^{m+1}]$ by
\eqref{eq:normal-moment}, followed by the sphere integral, gives
\begin{align*}
 &\Phi_3=-\frac{\operatorname{rk}(F)\,\pi(2m-1)!\operatorname{Vol}(S^{2m-1})}{2^{m}(m+1)!(m-2)!}\,\frac{m(2m+1)}{2(m+2)}h'(0)\\
 &\hspace{4em}\times [g(X,Y)Z_n-g(X,Z)Y_n+g(Y,Z)X_n]\,d\operatorname{vol}_{\partial M}.
\end{align*}

\paragraph{Case \(b\).}
Take \(r=-1,\ \ell=-2m+1,\ j=k=|\alpha|=0\). The coefficient is \(-i\).
\begin{align*}
\Phi_4={}&-i\int_{|\xi'|=1}\int_{\mathbb R}
 \operatorname{tr}_{S(TM)\otimes F}\!\Bigl[
 \pi^+_{\xi_n}\sigma_{-1}\!\left(c(X)c(Y)c(Z)\widetilde D_F^{-1}\right)\\
&\qquad{}\cdot \partial_{\xi_n}\sigma_{-2m+1}\!\left(\widetilde D_F^{-2m+2}\right)\Bigr](x_0,\xi)
 \,d\xi_n\,dS(\xi')\,d\operatorname{vol}_{\partial M}.
\end{align*}
Only $b$ and $c$ contain the perturbation. Integration by
parts moves $\partial_{\xi_n}$ to the projected left factor, changing
$-i$ into $+i$. The perturbation in the right symbol is
$-i(m-1)\{c(\xi),c(\omega)\}\otimes\Phi/(1+\xi_n^2)^m$,
where $\{U,V\}=UV+VU$. At fixed $\xi'$, the resulting scalar integrand is
\begin{equation}\label{co:btrace}
 -\frac{(m-1)\operatorname{tr}_F(\Phi)}
 {2(\xi_n-i)^2(1+\xi_n^2)^m}
 \operatorname{tr}_{S(TM)}
 \big[A(c(\xi')+i c(dx_n))\{c(\xi),c(\omega)\}\big].
\end{equation}
For clarity, the coefficient of a pair of directions $j,k$ in this
ordinary spinor trace follows directly by anticommuting the two indicated
generators. One obtains
\begin{align*}
 &2^{-m}\operatorname{tr}_{S(TM)}
       [A c(e_j)\{c(e_k),c(\omega)\}]\\
 ={}&-2\theta_1(e_k)[g(X,Y)Z_j-g(X,Z)Y_j+g(Y,Z)X_j]\\
 &-2[X_j\theta_3(e_k,Y,Z)-Y_j\theta_3(e_k,X,Z)
                          +Z_j\theta_3(e_k,X,Y)]\\
 &+2[g(X,Y)\theta_3(e_k,e_j,Z)-g(X,Z)\theta_3(e_k,e_j,Y)
                             +g(Y,Z)\theta_3(e_k,e_j,X)]\\
 &+2\theta_5(e_k,e_j,X,Y,Z).
\end{align*}
In particular the five-form can occur before angular integration. To
integrate \eqref{co:btrace}, the two required normal moments are
\begin{align*}
 -\frac12\int_{\mathbb R}\frac{d\xi_n}
    {(\xi_n-i)^2(1+\xi_n^2)^m}
 &=\frac{\pi(2m)!}{2^{2m+1}(m-1)!(m+1)!},\\
 -\frac12\int_{\mathbb R}\frac{\xi_n\,d\xi_n}
    {(\xi_n-i)^2(1+\xi_n^2)^m}
 &=-\frac{i\pi(2m)!}{2^{2m+1}m!(m+1)!}.
\end{align*}
The linear terms in $\xi_a$ integrate to zero and the angular average
of $\xi_a\xi_b$ is $\delta_{ab}/(2m)$. The surviving pairs have $j=k$,
which removes every $\theta_5$ term. The metric part and the part linear in the potential are
\begin{align*}
 &\left.\Phi_4\right|_{\omega=0}=\frac{\operatorname{rk}(F)\,\pi(2m-1)!\operatorname{Vol}(S^{2m-1})}{2^{m}(m+1)!(m-2)!}\,\frac{m(5m+7)}{2(m+2)}h'(0)\\
 &\hspace{4em}\times [g(X,Y)Z_n-g(X,Z)Y_n+g(Y,Z)X_n]\,d\operatorname{vol}_{\partial M}.
\end{align*}
\begin{align*}
 &\Phi_4-\left.\Phi_4\right|_{\omega=0}=\frac{\operatorname{tr}_F(\Phi)\,\pi(2m-1)!\operatorname{Vol}(S^{2m-1})}{2^{2m}(m+1)!(m-2)!}\\
 &\quad\times\operatorname{tr}_{S(TM)}\!\bigg[A\bigg(-(m+1)c(\omega)+\frac12\sum_{a=1}^{2m}c(e_a)c(\omega)c(e_a)\\
 &\hspace{8em}{}+c(dx_n)c(\omega)c(dx_n)\bigg)\bigg]\,d\operatorname{vol}_{\partial M}.
\end{align*}

\paragraph{Case \(c\).}
Take \(r=-2,\ \ell=-2m+2,\ j=k=|\alpha|=0\). The coefficient is \(-i\).
\begin{align*}
\Phi_5={}&-i\int_{|\xi'|=1}\int_{\mathbb R}
 \operatorname{tr}_{S(TM)\otimes F}\!\Bigl[
 \pi^+_{\xi_n}\sigma_{-2}\!\left(c(X)c(Y)c(Z)\widetilde D_F^{-1}\right)\\
&\qquad{}\cdot \partial_{\xi_n}\sigma_{-2m+2}\!\left(\widetilde D_F^{-2m+2}\right)\Bigr](x_0,\xi)
 \,d\xi_n\,dS(\xi')\,d\operatorname{vol}_{\partial M}.
\end{align*}
The part of the projected inverse symbol linear in
$c(\omega)\otimes\Phi$ is
\begin{align*}
 &-\frac{[c(\xi')+i c(dx_n)]c(\omega)[c(\xi')+i c(dx_n)]}
           {4(\xi_n-i)^2}\otimes\Phi\\
 &\qquad-\frac{i[c(\xi')c(\omega)c(\xi')
                    +c(dx_n)c(\omega)c(dx_n)]}{4(\xi_n-i)}
           \otimes\Phi .
\end{align*}
In the residue formula the insertion $A$ multiplies this symbol from the left, and the scalar
$-2(m-1)\xi_n/(1+\xi_n^2)^m$ multiplies it from the right. In the following identity
the factor $-i$ is the boundary coefficient of case $c$. One obtains
\[
 -i2^{m-1}(m-1)\int_{\mathbb R}
 \frac{\xi_n\,d\xi_n}{(1+\xi_n^2)^m(\xi_n-i)^q}
 =\frac{\pi(2m-1)!}{2^{m+1}(m+1)!(m-2)!}
 \begin{cases}
 -2i(m+1)/(2m-1), &q=1,\\
 2,&q=2.
 \end{cases}
\]
The normal moments and the identities for the first two angular moments now give
\begin{align*}
 &\left.\Phi_5\right|_{\omega=0}=-\frac{\operatorname{rk}(F)\,\pi(2m-1)!\operatorname{Vol}(S^{2m-1})}{2^{m}(m+1)!(m-2)!}\,\frac{m(2m+7)}{2(m+2)}h'(0)\\
 &\hspace{4em}\times [g(X,Y)Z_n-g(X,Z)Y_n+g(Y,Z)X_n]\,d\operatorname{vol}_{\partial M}.
\end{align*}
\begin{align*}
 &\Phi_5-\left.\Phi_5\right|_{\omega=0}=\frac{\operatorname{tr}_F(\Phi)\,\pi(2m-1)!\operatorname{Vol}(S^{2m-1})}{2^{2m}(m+1)!(m-2)!}\\
 &\quad\times\operatorname{tr}_{S(TM)}\!\bigg[A\bigg(\frac{3}{2(2m-1)}\sum_{a=1}^{2m}c(e_a)c(\omega)c(e_a)\\
 &\hspace{8em}{}+\frac{2-m}{2m-1}c(dx_n)c(\omega)c(dx_n)\bigg)\bigg]\,d\operatorname{vol}_{\partial M}.
\end{align*}
Adding the perturbation terms of cases $b$ and $c$ leaves, inside the
spinor trace, the Clifford factor
\[
 \frac{m+1}{2m-1}\left[\sum_{j=1}^{2m+1}c(e_j)c(\omega)c(e_j)
                            -(2m-1)c(\omega)\right].
\]
For an odd $q$-form $\beta$, anticommutation yields
$\sum_{j=1}^n c(e_j)c(\beta)c(e_j)=(n-2q)c(\beta)$.
This factor therefore annihilates degree one. Its multiplier on
degree three is $-4(m+1)/(2m-1)$. Higher odd degrees have zero
trace against $A$, since an even Clifford monomial in odd dimension
has zero ordinary trace unless it reduces to a scalar.
Finally $2^{-m}\operatorname{tr}_{S(TM)}[A c(\theta_3)]=\theta_3(X,Y,Z)$.
The metric terms are those of \eqref{dt:formula1} multiplied by
$\operatorname{rk}(F)$, with the normal derivative supplied by $\Phi_2$.
Integrating the boundary density proves \eqref{co:first}.
\end{proof}

\subsection{Boundary residue of \texorpdfstring{$\pi^+\!\bigl(c(X)c(Y)c(Z)\widetilde D_F^{-2}\bigr)\circ\pi^+\!\bigl(\widetilde D_F^{-2m+3}\bigr)$}{the Connes-type product with insertion c(X)c(Y)c(Z) and powers -2 and -2m+3}}

\begin{theorem}\label{co:two}
Under the same hypotheses,
\begin{equation}\label{co:second}
\begin{aligned}
 &\widetilde{\operatorname{Wres}}\!\left[\pi^+\!\left(c(X)c(Y)c(Z)\widetilde D_F^{-2}\right)\circ\pi^+\!\left(\widetilde D_F^{-2m+3}\right)\right]\\
 &\quad=\frac{3\pi(2m-2)!\operatorname{Vol}(S^{2m-1})}
              {2^m(m+1)!(m-2)!}\int_{\partial M}\!\bigg\{\\
 &\qquad\operatorname{rk}(F)
 \left.\partial_{x_n}\big[g(X,Y)Z_n-g(X,Z)Y_n+g(Y,Z)X_n\big]
 \right|_{x_n=0}\\
 &\quad-\frac m2\operatorname{rk}(F)h'(0)
          [g(X,Y)Z_n-g(X,Z)Y_n+g(Y,Z)X_n]\\
 &\quad-\frac{4(m+1)}3\operatorname{tr}_F(\Phi)
       (\omega_{[3]}+\kappa_n*\omega_{[n-3]})(X,Y,Z)\bigg\}
             \,d\operatorname{vol}_{\partial M}.
\end{aligned}
\end{equation}
\end{theorem}

\begin{proof}
As in the first residue, the geometric contribution is the unperturbed
spinor density multiplied by $\operatorname{rk}(F)$. The potential occurs
in the two subleading symbols.
\paragraph{Case \(a_1\).}
Take \(r=-2,\ \ell=-2m+3,\ j=k=0,\ |\alpha|=1\). The coefficient is \(-1\).
\begin{align*}
\widetilde\Phi_1={}&-\int_{|\xi'|=1}\int_{\mathbb R}\sum_{a=1}^{2m}
 \operatorname{tr}_{S(TM)\otimes F}\!\Bigl[
 \partial_{\xi_a}\pi^+_{\xi_n}\sigma_{-2}\!\left(c(X)c(Y)c(Z)\widetilde D_F^{-2}\right)\\
&\qquad{}\cdot \partial_{x_a}\partial_{\xi_n}\sigma_{-2m+3}\!\left(\widetilde D_F^{-2m+3}\right)\Bigr](x_0,\xi)
 \,d\xi_n\,dS(\xi')\,d\operatorname{vol}_{\partial M}.
\end{align*}
The right principal symbol is $i c(\xi)(|\xi|^2)^{1-m}$.
In the frame over tangential geodesic coordinates,
$\partial_{x_a}c(\xi)(x_0)=0$ and $\partial_{x_a}|\xi|^2(x_0)=0$
for $a\leq2m$. The mixed derivative on the right vanishes, so $\widetilde\Phi_1=0$.

\paragraph{Case \(a_2\).}
Take \(r=-2,\ \ell=-2m+3,\ j=1,\ k=|\alpha|=0\). The coefficient is \(-\frac12\).
\begin{align*}
\widetilde\Phi_2={}&-\frac12\int_{|\xi'|=1}\int_{\mathbb R}
 \operatorname{tr}_{S(TM)\otimes F}\!\Bigl[
 \partial_{x_n}\pi^+_{\xi_n}\sigma_{-2}\!\left(c(X)c(Y)c(Z)\widetilde D_F^{-2}\right)\\
&\qquad{}\cdot \partial_{\xi_n}^{2}\sigma_{-2m+3}\!\left(\widetilde D_F^{-2m+3}\right)\Bigr](x_0,\xi)
 \,d\xi_n\,dS(\xi')\,d\operatorname{vol}_{\partial M}.
\end{align*}
Here the projected scalar inverse symbol has both a moving pole
and a variable prefactor. Differentiating before boundary restriction gives
\[
 -\frac{i\,\partial_{x_n}[c(X)c(Y)c(Z)]}{2(\xi_n-i)}
 +h'(0)A\left[\frac{i}{4(\xi_n-i)}+\frac1{4(\xi_n-i)^2}\right].
\]
In the right symbol the coefficient of $c(dx_n)$ is
$i\partial_{\xi_n}^2[\xi_n(1+\xi_n^2)^{1-m}]$.
The remaining Clifford component is linear in $\xi'$ and integrates
to zero on the sphere. Two integrations by parts in $\xi_n$ and
\eqref{eq:normal-moment} give
\begin{align*}
 &\widetilde\Phi_2=\frac{3\operatorname{rk}(F)\,\pi(2m-2)!\operatorname{Vol}(S^{2m-1})}{2^{m}(m+1)!(m-2)!}\bigg\{\\
 &\hspace{4em}\left.\partial_{x_n}\!\left[g(X,Y)Z_n-g(X,Z)Y_n+g(Y,Z)X_n\right]\right|_{x_n=0}\\
 &\hspace{4em}-\frac{5m}{2(m+2)}h'(0)[g(X,Y)Z_n-g(X,Z)Y_n+g(Y,Z)X_n]\bigg\}\,d\operatorname{vol}_{\partial M}.
\end{align*}

\paragraph{Case \(a_3\).}
Take \(r=-2,\ \ell=-2m+3,\ j=|\alpha|=0,\ k=1\). The coefficient is \(-\frac12\).
\begin{align*}
\widetilde\Phi_3={}&-\frac12\int_{|\xi'|=1}\int_{\mathbb R}
 \operatorname{tr}_{S(TM)\otimes F}\!\Bigl[
 \partial_{\xi_n}\pi^+_{\xi_n}\sigma_{-2}\!\left(c(X)c(Y)c(Z)\widetilde D_F^{-2}\right)\\
&\qquad{}\cdot \partial_{\xi_n}\partial_{x_n}\sigma_{-2m+3}\!\left(\widetilde D_F^{-2m+3}\right)\Bigr](x_0,\xi)
 \,d\xi_n\,dS(\xi')\,d\operatorname{vol}_{\partial M}.
\end{align*}
The differentiated left projection is $i A/[2(\xi_n-i)^2]$.
After differentiating the right symbol in the collar and restricting to
$x_n=0$ and $|\xi'|_{g^{\partial M}}=1$, one obtains
\[
 \frac{i(m-1)h'(0)}{(1+\xi_n^2)^{m+1}}
 \big[\xi_n(2m-1-\xi_n^2)c(\xi')
       +((2m-1)\xi_n^2-1)c(dx_n)\big].
\]
The component with $c(\xi')$ disappears on angular integration.
The normal component gives the density
\begin{align*}
 &\widetilde\Phi_3=\frac{3\operatorname{rk}(F)\,\pi(2m-2)!\operatorname{Vol}(S^{2m-1})}{2^{m}(m+1)!(m-2)!}\,\frac{m(2m-1)}{2(m+2)}h'(0)\\
 &\hspace{4em}\times [g(X,Y)Z_n-g(X,Z)Y_n+g(Y,Z)X_n]\,d\operatorname{vol}_{\partial M}.
\end{align*}

\paragraph{Case \(b\).}
Take \(r=-3,\ \ell=-2m+3,\ j=k=|\alpha|=0\). The coefficient is \(-i\).
\begin{align*}
\widetilde\Phi_4={}&-i\int_{|\xi'|=1}\int_{\mathbb R}
 \operatorname{tr}_{S(TM)\otimes F}\!\Bigl[
 \pi^+_{\xi_n}\sigma_{-3}\!\left(c(X)c(Y)c(Z)\widetilde D_F^{-2}\right)\\
&\qquad{}\cdot \partial_{\xi_n}\sigma_{-2m+3}\!\left(\widetilde D_F^{-2m+3}\right)\Bigr](x_0,\xi)
 \,d\xi_n\,dS(\xi')\,d\operatorname{vol}_{\partial M}.
\end{align*}
The left subleading symbol contains
$-iA\{c(\xi),c(\omega)\}\otimes\Phi/(1+\xi_n^2)^2$.
The projection acts separately on each $\xi_j/(1+\xi_n^2)^2$. In
particular
\[
 \pi^+_{\xi_n}\frac{1}{(1+\xi_n^2)^2}
 =-\frac{1}{4(\xi_n-i)^2}-\frac{i}{4(\xi_n-i)},\qquad
 \pi^+_{\xi_n}\frac{\xi_n}{(1+\xi_n^2)^2}
 =-\frac{i}{4(\xi_n-i)^2}.
\]
The perturbation in the trace at fixed $\xi'$ including the boundary
coefficient is
\begin{align*}
 &-i\operatorname{tr}_F(\Phi)\sum_{j,k=1}^n
 \pi^+_{\xi_n}\!\left(\frac{\xi_j}{(1+\xi_n^2)^2}\right)
 \partial_{\xi_n}\!\left(\frac{\xi_k}{(1+\xi_n^2)^{m-1}}\right)
 \operatorname{tr}_{S(TM)}[A\{c(e_j),c(\omega)\}c(e_k)] .
\end{align*}
For the tangential and normal diagonal pairs the two normal integrals
are, respectively,
\begin{align*}
 -i\int_{\mathbb R}\!
 \left(\pi^+_{\xi_n}\frac{1}{(1+\xi_n^2)^2}\right)
 \partial_{\xi_n}(1+\xi_n^2)^{1-m}\,d\xi_n
 &=\frac{3m\pi(2m-2)!}{2^{2m}(m-2)!(m+1)!},\\
 -i\int_{\mathbb R}\!
 \left(\pi^+_{\xi_n}\frac{\xi_n}{(1+\xi_n^2)^2}\right)
 \partial_{\xi_n}\big[\xi_n(1+\xi_n^2)^{1-m}\big]\,d\xi_n
 &=\frac{3\pi(2m-2)!}{2^{2m}(m-2)!(m+1)!}.
\end{align*}
The angular integration multiplies a tangential diagonal term by
$\operatorname{Vol}(S^{2m-1})/(2m)$ and a normal diagonal term by
$\operatorname{Vol}(S^{2m-1})$. The surviving matrix in direction $j$
is $c(e_j)c(\omega)c(e_j)-c(\omega)$. The two contributions are
\begin{align*}
 &\left.\widetilde\Phi_4\right|_{\omega=0}=\frac{3\operatorname{rk}(F)\,\pi(2m-2)!\operatorname{Vol}(S^{2m-1})}{2^{m}(m+1)!(m-2)!}\,\frac{m(m+7)}{2(m+2)}h'(0)\\
 &\hspace{4em}\times [g(X,Y)Z_n-g(X,Z)Y_n+g(Y,Z)X_n]\,d\operatorname{vol}_{\partial M}.
\end{align*}
\begin{align*}
 &\widetilde\Phi_4-\left.\widetilde\Phi_4\right|_{\omega=0}=\frac{3\operatorname{tr}_F(\Phi)\,\pi(2m-2)!\operatorname{Vol}(S^{2m-1})}{2^{2m}(m+1)!(m-2)!}\\
 &\quad\times\operatorname{tr}_{S(TM)}\!\bigg[A\bigg(-(m+1)c(\omega)+\frac12\sum_{a=1}^{2m}c(e_a)c(\omega)c(e_a)\\
 &\hspace{8em}{}+c(dx_n)c(\omega)c(dx_n)\bigg)\bigg]\,d\operatorname{vol}_{\partial M}.
\end{align*}

\paragraph{Case \(c\).}
Take \(r=-2,\ \ell=-2m+2,\ j=k=|\alpha|=0\). The coefficient is \(-i\).
\begin{align*}
\widetilde\Phi_5={}&-i\int_{|\xi'|=1}\int_{\mathbb R}
 \operatorname{tr}_{S(TM)\otimes F}\!\Bigl[
 \pi^+_{\xi_n}\sigma_{-2}\!\left(c(X)c(Y)c(Z)\widetilde D_F^{-2}\right)\\
&\qquad{}\cdot \partial_{\xi_n}\sigma_{-2m+2}\!\left(\widetilde D_F^{-2m+3}\right)\Bigr](x_0,\xi)
 \,d\xi_n\,dS(\xi')\,d\operatorname{vol}_{\partial M}.
\end{align*}
Write
$\widetilde D_F^{-2m+3}=(\widetilde D_F^2)^{1-m}\widetilde D_F$.
Integration by parts changes the boundary coefficient to $+i$, and
the differentiated left projection is $iA/[2(\xi_n-i)^2]$.
At fixed $\xi'$ the scalar integrand contributed by the perturbation is
\[
 -\frac{\operatorname{tr}_F(\Phi)}{2(\xi_n-i)^2}
 \operatorname{tr}_{S(TM)}\!\left[
 A\left(\frac{c(\omega)}{(1+\xi_n^2)^{m-1}}
 +\frac{(m-1)\{c(\xi),c(\omega)\}c(\xi)}{(1+\xi_n^2)^m}\right)\right].
\]
The two moments in the preceding proof apply with $m$ replaced by $m-1$
or $m$. The additional normal diagonal moment is
\[
 -\frac12\int_{\mathbb R}
 \frac{\xi_n^2\,d\xi_n}{(\xi_n-i)^2(1+\xi_n^2)^m}
 =\frac{\pi(2m-2)!(m-2)}{2^{2m}(m-1)!(m+1)!}.
\]
Combining these moments with the angular averages gives
\begin{align*}
 &\left.\widetilde\Phi_5\right|_{\omega=0}=-\frac{3\operatorname{rk}(F)\,\pi(2m-2)!\operatorname{Vol}(S^{2m-1})}{2^{m}(m+1)!(m-2)!}\,\frac{m(4m+3)}{2(m+2)}h'(0)\\
 &\hspace{4em}\times [g(X,Y)Z_n-g(X,Z)Y_n+g(Y,Z)X_n]\,d\operatorname{vol}_{\partial M}.
\end{align*}
\begin{align*}
 &\widetilde\Phi_5-\left.\widetilde\Phi_5\right|_{\omega=0}=\frac{3\operatorname{tr}_F(\Phi)\,\pi(2m-2)!\operatorname{Vol}(S^{2m-1})}{2^{2m}(m+1)!(m-2)!}\\
 &\quad\times\operatorname{tr}_{S(TM)}\!\bigg[A\bigg(-\frac{2(m-2)(m+1)}3c(\omega)\\
 &\hspace{8em}{}+\frac{2m-1}{6}\sum_{a=1}^{2m}c(e_a)c(\omega)c(e_a)\\
 &\hspace{8em}{}+\frac{m-2}{3}c(dx_n)c(\omega)c(dx_n)\bigg)\bigg]\,d\operatorname{vol}_{\partial M}.
\end{align*}
Adding the terms linear in the potential from cases $b$ and $c$ gives
the Clifford factor
\[
 \frac{m+1}{3}\left[\sum_{j=1}^{2m+1}c(e_j)c(\omega)c(e_j)
                                -(2m-1)c(\omega)\right].
\]
Its eigenvalues on one- and three-forms are $0$ and $-4(m+1)/3$,
respectively. Adding the metric terms, with the normal derivative
supplied by $\widetilde\Phi_2$, and integrating over $\partial M$
gives \eqref{co:second}.
\end{proof}

\subsection{The three-form contribution of the Clifford potential}
The scalar coefficients in \eqref{co:first} and \eqref{co:second} are defined on the entire boundary. Multiplication by $d\operatorname{vol}_{\partial M}$ gives the respective residue densities.
The twist enters the metric term through its fiber dimension and the
Clifford potential through its fiber trace. Thus a pointwise trace-free
$\Phi$ has the same two residues as the unperturbed twisted Dirac
operator. For the potential $c(\omega)=c(T)+i c(Y_0)$ of
\cite[Example~2.6]{WW}, the three-form in \eqref{co:dual} is exactly $T$.
Both residues are insensitive to the vector field $Y_0$. For a general
potential they measure the combination
$\operatorname{tr}_F(\Phi)(\omega_{[3]}+\kappa_n*\omega_{[n-3]})$.
These formulas give a boundary response linear in the three-form
coupling of a twisted Dirac field. The curvature of the connection on
$F$ does not enter at this symbol order.

For arbitrary $X,Y,Z$, the linear combination that cancels the metric
contribution is
\begin{align}
 &\frac{3}{2m+2}\widetilde{\operatorname{Wres}}\!\left[\pi^+\!\left(A\widetilde D_F^{-1}\right)\circ\pi^+\!\left(\widetilde D_F^{-2m+2}\right)\right]\notag\\
 &\quad{}+\frac{2m-1}{2m+2}\widetilde{\operatorname{Wres}}\!\left[\pi^+\!\left(A\widetilde D_F^{-2}\right)\circ\pi^+\!\left(\widetilde D_F^{-2m+3}\right)\right]
 \notag\\
 &\quad=-\frac{4\pi(2m-2)!\operatorname{Vol}(S^{2m-1})}
                   {2^m m!(m-2)!}
 \int_{\partial M}\operatorname{tr}_F(\Phi)
 (\omega_{[3]}+\kappa_n*\omega_{[n-3]})(X,Y,Z)
 \,d\operatorname{vol}_{\partial M}.\label{co:weighted}
\end{align}
Here integration pairs a boundary scalar with the induced volume
density. The evaluated ambient three-form may have a normal slot when
the fields have normal components. The difference of the two residues
instead cancels the three-form pointwise. If the contraction
$g(X,Y)Z_n-g(X,Z)Y_n+g(Y,Z)X_n$ equals a known constant $C_0$ throughout
a collar, the inward trace mean curvature $K=-mh'(0)$ gives
\begin{align*}
 &\widetilde{\operatorname{Wres}}\!\left[\pi^+\!\left(A\widetilde D_F^{-2}\right)\circ\pi^+\!\left(\widetilde D_F^{-2m+3}\right)\right]\\
 &\quad{}-\widetilde{\operatorname{Wres}}\!\left[\pi^+\!\left(A\widetilde D_F^{-1}\right)\circ\pi^+\!\left(\widetilde D_F^{-2m+2}\right)\right]\\
 &\qquad=\frac{\operatorname{rk}(F)\pi(2m-2)!
                   \operatorname{Vol}(S^{2m-1})}{2^m m!(m-2)!}\,
 C_0\int_{\partial M}K\,d\operatorname{vol}_{\partial M}.
\end{align*}
For example $X=Y=Z=\partial_{x_n}$ on the collar gives $C_0=1$.
The difference then determines the total mean curvature with the inward
normal convention. Prescribing the contraction only on the boundary
leaves its normal derivative free.

\begin{corollary}
If $X,Y,Z$ are tangent on $\partial M$, the term proportional to
$h'(0)$ vanishes. The derivative term in both theorems reduces to
\[
 \left.[g(X,Y)\partial_{x_n}Z_n-g(X,Z)\partial_{x_n}Y_n
                  +g(Y,Z)\partial_{x_n}X_n]\right|_{x_n=0}.
\]
If the extensions are tangent to every collar slice, this derivative
vanishes as well. In that case both residues equal the right-hand side
of \eqref{co:weighted}.
\end{corollary}

\paragraph{Clifford duality and the information retained by the observable.}
Set
\[
 \Theta=\operatorname{tr}_F(\Phi)
       (\omega_{[3]}+\kappa_n*\omega_{[n-3]}).
\]
Localizing \eqref{co:weighted} in one insertion determines $\Theta$ as an ambient three-form along the boundary, including components with one normal argument. It does not determine its two summands separately. For example, for any $(n-3)$-form $\eta$, the simultaneous changes
\[
 \omega_{[n-3]}\longmapsto\omega_{[n-3]}+\eta,
 \qquad \omega_{[3]}\longmapsto\omega_{[3]}-\kappa_n*\eta
\]
leave $\Theta$ unchanged. The two degrees are distinct for odd $n$. This ambiguity is the direct algebraic consequence of the scalar action of the odd Clifford volume element, rather than an additional dynamical duality assumption. In dimension five the selected combination is $\omega_{[3]}+i*\omega_{[2]}$; in dimension seven it is $\omega_{[3]}+*\omega_{[4]}$, with the convention fixed in Section~\ref{sec:prelim}.

\paragraph{Internal trace and gauge invariance.}
Under a unitary change of frame $U$ in $F$, the endomorphism transforms to $U\Phi U^{-1}$ and its trace is unchanged. The potential contribution is consequently independent of this choice of frame. Writing $r=\operatorname{rk}(F)>0$ and
\[
 \Phi=\frac{\operatorname{tr}_F\Phi}{r}\operatorname{Id}_F+\Phi_0,
 \qquad\operatorname{tr}_F\Phi_0=0,
\]
shows that the observable retains only the scalar part in the internal factor. This is a selection rule of the ordinary fiber trace. A trace-free interaction can still be nonzero and affect other spectral invariants. Nor should the absence of the twisting curvature here be interpreted as absence of gauge dynamics: the calculation probes only the indicated symbol order, whereas spectral actions involve other coefficients \cite{CC,Vass}.

\paragraph{Boundary response versus a physical action.}
The weighted functional can be interpreted as the boundary response selected by three Clifford probes of a form-valued interaction. Its dependence on $\Phi$ and on the effective three-form is explicitly linear in the zero-order potential at this order. The difference of the two residues independently supplies the metric response, multiplied by the number of internal degrees of freedom. For general complex $\omega$, $\Theta$ and the residue may be complex. A real physical action would require appropriate reality and adjointness conditions and, in general, a specified boundary realization. The present construction establishes a geometric observable of the differential expression without imposing these further physical choices.

\section{The de Rham--Hodge operator}\label{sec:hodge}

The two residues for \(D_{\mathrm{dR}}=d+\delta\) are defined on the full
complex exterior bundle. This realization exists on an oriented manifold
without choosing a spin structure. The exterior connection contributes
additional terms to the local symbols. We identify these terms before
evaluating the ordinary trace.

\subsection{Boundary residue of \texorpdfstring{$\pi^+\!\bigl(c(X)c(Y)c(Z)D_{\mathrm{dR}}^{-1}\bigr)\circ\pi^+\!\bigl(D_{\mathrm{dR}}^{-2m+2}\bigr)$}{the de Rham--Hodge product with insertion c(X)c(Y)c(Z) and powers -1 and -2m+2}}\label{hd:firstsub}

Let \(\widehat c(U)\) denote the sum of exterior multiplication by
\(g(U,\cdot)\) and contraction by \(U\). The operator \(c(U)\) is the
difference of the same two operations \cite[Sec.~2]{LWW}. Thus
\[
 c(U)\widehat c(W)+\widehat c(W)c(U)=0,\qquad
 \widehat c(U)\widehat c(W)+\widehat c(W)\widehat c(U)=2g(U,W).
\]
At \(x_0\), the matrix of the induced connection in the direction \(e_a\)
in the normally parallel exterior frame is
\[
 \frac{h'(0)}4\big[c(e_a)c(dx_n)-\widehat c(e_a)\widehat c(dx_n)\big],
 \qquad 1\le a\le 2m.
\]
The normal connection matrix vanishes. Multiplying these connection
matrices by \(c(e_a)\) and summing, and then using the symbol product for
the square, gives
\begin{align}
 \sigma_0(D_{\mathrm{dR}})(x_0)
 &=-\frac m2 h'(0)c(dx_n)
 -\frac{h'(0)}4\sum_{a=1}^{2m}c(e_a)\widehat c(e_a)\widehat c(dx_n),
 \label{hd:zero}\\
 \sigma_1(D_{\mathrm{dR}}^2)(x_0,\xi)
 &=i h'(0)\left[m\xi_n+\frac12 c(dx_n)c(\xi')
 +\frac12\sum_{a=1}^{2m}\xi_a\widehat c(e_a)\widehat c(dx_n)\right].
 \label{hd:square}
\end{align}
In \eqref{hd:square}, the middle term comes from
\(i c(dx_n)\partial_{x_n}c(\xi)\).

Each term with a hatted pair in either subleading symbol enters the
boundary formula exactly once. All other factors are ordinary Clifford
matrices, and their total number, including the insertion, is even.
Conjugation by \(\widehat c(dx_n)\) changes the sign of
\(\widehat c(e_a)\widehat c(dx_n)\) and preserves an even product of
ordinary Clifford factors. Every such mixed term therefore has zero
ordinary trace at fixed \(\xi'\) and \(\xi_n\), including after projection.
The same argument applies to both factorizations. The remaining
Clifford reductions agree, term by term, with the untwisted spin
calculation in Section~\ref{sec:torsion}, with ordinary trace evaluated
on a space of dimension \(2^{2m+1}\). In particular,
\begin{equation}\label{hd:trace}
 \operatorname{tr}_{\Lambda}[c(X)c(Y)c(Z)c(W)]
 =2^{2m+1}\big[g(X,Y)g(Z,W)-g(X,Z)g(Y,W)+g(Y,Z)g(X,W)\big].
\end{equation}

For the factors of orders \(-1\) and \(-2m+2\), the leading symbols are
\(i c(X)c(Y)c(Z)c(\xi)/|\xi|^2\) and \((|\xi|^2)^{1-m}\).
\paragraph{Case \(a_1\).}
Take \(r=-1,\ \ell=-2m+2,\ j=k=0,\ |\alpha|=1\). The coefficient is \(-1\).
\begin{align*}
\Phi_1={}&-\int_{|\xi'|=1}\int_{\mathbb R}\sum_{a=1}^{2m}
 \operatorname{tr}_{\Lambda}\!\Bigl[
 \partial_{\xi_a}\pi^+_{\xi_n}\sigma_{-1}\!\left(c(X)c(Y)c(Z)D_{\mathrm{dR}}^{-1}\right)\\
&\qquad{}\cdot \partial_{x_a}\partial_{\xi_n}\sigma_{-2m+2}\!\left(D_{\mathrm{dR}}^{-2m+2}\right)\Bigr](x_0,\xi)
 \,d\xi_n\,dS(\xi')\,d\operatorname{vol}_{\partial M}.
\end{align*}
For every tangential coordinate index \(a\), the right factor satisfies
\[
 \partial_{x_a}\partial_{\xi_n}(|\xi|^2)^{1-m}
 =2m(m-1)\xi_n(|\xi|^2)^{-m-1}\partial_{x_a}|\xi|^2=0
 \quad\text{at }x_0 .
\]
Thus the integrand vanishes for every normal covariable.
Both the normal integral at fixed \(\xi'\) and the resulting boundary
density \(\Phi_1\) are zero.

\paragraph{Case \(a_2\).}
Take \(r=-1,\ \ell=-2m+2,\ j=1,\ k=|\alpha|=0\). The coefficient is \(-\frac12\).
\begin{align*}
\Phi_2={}&-\frac12\int_{|\xi'|=1}\int_{\mathbb R}
 \operatorname{tr}_{\Lambda}\!\Bigl[
 \partial_{x_n}\pi^+_{\xi_n}\sigma_{-1}\!\left(c(X)c(Y)c(Z)D_{\mathrm{dR}}^{-1}\right)\\
&\qquad{}\cdot \partial_{\xi_n}^{2}\sigma_{-2m+2}\!\left(D_{\mathrm{dR}}^{-2m+2}\right)\Bigr](x_0,\xi)
 \,d\xi_n\,dS(\xi')\,d\operatorname{vol}_{\partial M}.
\end{align*}
The left projection throughout the collar is
\[
 \frac{c(X)c(Y)c(Z)}{2(\xi_n-i|\xi'|_g)}
 \left[\frac{c(\xi')}{|\xi'|_g}+i c(dx_n)\right].
\]
In the chosen frame the quotient inside the brackets is constant along
the normal, with the coordinate covariable fixed. Differentiating both
the insertion and the pole, then restricting to \(x_0,\ |\xi'|=1\), gives
\begin{align*}
 &\frac{\partial_{x_n}[c(X)c(Y)c(Z)]\,[c(\xi')+i c(dx_n)]}{2(\xi_n-i)}
 +\frac{h'(0)c(X)c(Y)c(Z)[i c(\xi')-c(dx_n)]}{4(\xi_n-i)^2}.
\end{align*}
The second normal covariable derivative of the right symbol is
\(2(m-1)[(2m-1)\xi_n^2-1]/(1+\xi_n^2)^{m+1}\).
Pairing the two derivatives uses the trace in \eqref{hd:trace}.
Multiplying by \(-1/2\) and using \eqref{eq:normal-moment} gives the
following integral over \(\xi_n\) at fixed unit \(\xi'\):
\begin{align*}
 &\frac{\pi(2m-1)!}{(m+1)!(m-2)!}\bigg\{
 -2\partial_{x_n}[g(X,Y)Z_n-g(X,Z)Y_n+g(Y,Z)X_n]\\
 &\qquad{}+2i\sum_{a=1}^{2m}\xi_a\,\partial_{x_n}[g(X,Y)Z_a-g(X,Z)Y_a+g(Y,Z)X_a]\\
 &\qquad{}+\frac{3m h'(0)}{m+2}[g(X,Y)Z_n-g(X,Z)Y_n+g(Y,Z)X_n]\\
 &\qquad{}-\frac{3im h'(0)}{m+2}\sum_{a=1}^{2m}\xi_a[g(X,Y)Z_a-g(X,Z)Y_a+g(Y,Z)X_a]\bigg\}.
\end{align*}
All components in the integrated expressions are evaluated at \(x_0\).
The frame is parallel along the normal and the tangential coefficients
\(\xi_a\) appearing after differentiation are fixed. Upon angular
integration the odd sums vanish, so
\begin{align*}
 \Phi_2={}&\frac{\pi(2m-1)!\operatorname{Vol}(S^{2m-1})}{(m+1)!(m-2)!}\bigg\{
 -2\partial_{x_n}[g(X,Y)Z_n-g(X,Z)Y_n+g(Y,Z)X_n]\\
 &\qquad{}+\frac{3m h'(0)}{m+2}[g(X,Y)Z_n-g(X,Z)Y_n+g(Y,Z)X_n]\bigg\}\,d\operatorname{vol}_{\partial M}.
\end{align*}

\paragraph{Case \(a_3\).}
Take \(r=-1,\ \ell=-2m+2,\ j=|\alpha|=0,\ k=1\). The coefficient is \(-\frac12\).
\begin{align*}
\Phi_3={}&-\frac12\int_{|\xi'|=1}\int_{\mathbb R}
 \operatorname{tr}_{\Lambda}\!\Bigl[
 \partial_{\xi_n}\pi^+_{\xi_n}\sigma_{-1}\!\left(c(X)c(Y)c(Z)D_{\mathrm{dR}}^{-1}\right)\\
&\qquad{}\cdot \partial_{\xi_n}\partial_{x_n}\sigma_{-2m+2}\!\left(D_{\mathrm{dR}}^{-2m+2}\right)\Bigr](x_0,\xi)
 \,d\xi_n\,dS(\xi')\,d\operatorname{vol}_{\partial M}.
\end{align*}
The derivative of the left projected principal symbol is
\[
 -\frac{c(X)c(Y)c(Z)[c(\xi')+i c(dx_n)]}{2(\xi_n-i)^2}.
\]
Differentiating \((|\xi|^2)^{1-m}\) first along
\(x_n\) and then along \(\xi_n\) gives
\(2m(m-1)h'(0)\xi_n/(1+\xi_n^2)^{m+1}\).
The four-factor trace in \eqref{hd:trace}, the coefficient \(-1/2\),
and the normal moment \eqref{eq:normal-moment} give, at fixed \(\xi'\),
\begin{align*}
 &-\frac{\pi(2m-1)!}{(m+1)!(m-2)!}\,\frac{m(2m+1)}{m+2}\,h'(0)\bigg\{[g(X,Y)Z_n-g(X,Z)Y_n+g(Y,Z)X_n]\\
 &\hspace{41mm}{}-i\sum_{a=1}^{2m}\xi_a[g(X,Y)Z_a-g(X,Z)Y_a+g(Y,Z)X_a]\bigg\}.
\end{align*}
Integrating the tangentially constant part gives
\begin{align*}
 \Phi_3={}&-\frac{\pi(2m-1)!\operatorname{Vol}(S^{2m-1})}{(m+1)!(m-2)!}\,\frac{m(2m+1)}{m+2}\,h'(0)\\
 &\quad{}\times [g(X,Y)Z_n-g(X,Z)Y_n+g(Y,Z)X_n]\,d\operatorname{vol}_{\partial M}.
\end{align*}

\paragraph{Case \(b\).}
Take \(r=-1,\ \ell=-2m+1,\ j=k=|\alpha|=0\). The coefficient is \(-i\).
\begin{align*}
\Phi_4={}&-i\int_{|\xi'|=1}\int_{\mathbb R}
 \operatorname{tr}_{\Lambda}\!\Bigl[
 \pi^+_{\xi_n}\sigma_{-1}\!\left(c(X)c(Y)c(Z)D_{\mathrm{dR}}^{-1}\right)\\
&\qquad{}\cdot \partial_{\xi_n}\sigma_{-2m+1}\!\left(D_{\mathrm{dR}}^{-2m+2}\right)\Bigr](x_0,\xi)
 \,d\xi_n\,dS(\xi')\,d\operatorname{vol}_{\partial M}.
\end{align*}
The right subleading symbol is
\begin{align*}
 -(m-1)(|\xi|^2)^{-m}\sigma_1(D_{\mathrm{dR}}^2)
 -\frac{im(m-1)}{2}(|\xi|^2)^{-m-1}
    \sum_{j=1}^n\partial_{\xi_j}|\xi|^2\,\partial_{x_j}|\xi|^2 .
\end{align*}
Integration by parts in \(\xi_n\) transfers the derivative to the left
projection and changes the prefactor from \(-i\) to \(+i\). At fixed
unit \(\xi'\), the resulting scalar integrand, including this coefficient,
is
\begin{align*}
 &-\frac{2^{2m}(m-1)h'(0)}{(\xi_n-i)^2(1+\xi_n^2)^m}
 \left(m\xi_n-\frac i2+\frac{m\xi_n}{1+\xi_n^2}\right)\\
 &\quad{}\times\bigg\{
 \sum_{a=1}^{2m}\xi_a[g(X,Y)Z_a-g(X,Z)Y_a+g(Y,Z)X_a]
 +i[g(X,Y)Z_n-g(X,Z)Y_n+g(Y,Z)X_n]\bigg\}.
\end{align*}
Integrating this expression in \(\xi_n\) by \eqref{eq:normal-moment} yields
\begin{align*}
 &\frac{\pi(2m-1)!}{(m+1)!(m-2)!}\,\frac{m(5m+7)}{m+2}\,h'(0)\bigg\{[g(X,Y)Z_n-g(X,Z)Y_n+g(Y,Z)X_n]\\
 &\hspace{41mm}{}-i\sum_{a=1}^{2m}\xi_a[g(X,Y)Z_a-g(X,Z)Y_a+g(Y,Z)X_a]\bigg\}.
\end{align*}
The part linear in the tangential covariable integrates to zero, leaving
\begin{align*}
 \Phi_4={}&\frac{\pi(2m-1)!\operatorname{Vol}(S^{2m-1})}{(m+1)!(m-2)!}\,\frac{m(5m+7)}{m+2}\,h'(0)\\
 &\quad{}\times [g(X,Y)Z_n-g(X,Z)Y_n+g(Y,Z)X_n]\,d\operatorname{vol}_{\partial M}.
\end{align*}

\paragraph{Case \(c\).}
Take \(r=-2,\ \ell=-2m+2,\ j=k=|\alpha|=0\). The coefficient is \(-i\).
\begin{align*}
\Phi_5={}&-i\int_{|\xi'|=1}\int_{\mathbb R}
 \operatorname{tr}_{\Lambda}\!\Bigl[
 \pi^+_{\xi_n}\sigma_{-2}\!\left(c(X)c(Y)c(Z)D_{\mathrm{dR}}^{-1}\right)\\
&\qquad{}\cdot \partial_{\xi_n}\sigma_{-2m+2}\!\left(D_{\mathrm{dR}}^{-2m+2}\right)\Bigr](x_0,\xi)
 \,d\xi_n\,dS(\xi')\,d\operatorname{vol}_{\partial M}.
\end{align*}
The two terms of
\(\pi^+_{\xi_n}\sigma_{-2}(D_{\mathrm{dR}}^{-1})\) that survive
multiplication by the insertion and passage to the trace are
\[
 -\frac{h'(0)[c(\xi')+i c(dx_n)]}{4(\xi_n-i)^3}
 +\frac{m h'(0)[c(dx_n)-i c(\xi')]}{4(\xi_n-i)^2}.
\]
They multiply \(-2(m-1)\xi_n/(1+\xi_n^2)^m\) on the right, with the
boundary prefactor \(-i\).
Substitute these factors into the trace and apply \eqref{eq:normal-moment}
to the double and triple poles. Including \(-i\), the normal integral is
\begin{align*}
 &-\frac{\pi(2m-1)!}{(m+1)!(m-2)!}\,\frac{m(2m+7)}{m+2}\,h'(0)\bigg\{[g(X,Y)Z_n-g(X,Z)Y_n+g(Y,Z)X_n]\\
 &\hspace{41mm}{}-i\sum_{a=1}^{2m}\xi_a[g(X,Y)Z_a-g(X,Z)Y_a+g(Y,Z)X_a]\bigg\}.
\end{align*}
Angular integration gives
\begin{align*}
 \Phi_5={}&-\frac{\pi(2m-1)!\operatorname{Vol}(S^{2m-1})}{(m+1)!(m-2)!}\,\frac{m(2m+7)}{m+2}\,h'(0)\\
 &\quad{}\times [g(X,Y)Z_n-g(X,Z)Y_n+g(Y,Z)X_n]\,d\operatorname{vol}_{\partial M}.
\end{align*}

Summing the five cases gives the following density. Its normal
derivative term comes entirely from \(\Phi_2\).

\begin{theorem}\label{hd:first}
The boundary density of the first de Rham--Hodge residue is
\begin{align*}
 &-\frac{2\pi(2m-1)!\operatorname{Vol}(S^{2m-1})}{(m+1)!(m-2)!}
 \bigg\{\left.\partial_{x_n}[g(X,Y)Z_n-g(X,Z)Y_n+g(Y,Z)X_n]\right|_{\partial M}\\
 &\hspace{29mm}{}-\frac m2 h'(0)
 [g(X,Y)Z_n-g(X,Z)Y_n+g(Y,Z)X_n]\bigg\}\,d\operatorname{vol}_{\partial M}.
\end{align*}
The integral of this density over \(\partial M\) equals

\[
 \widetilde{\operatorname{Wres}}\!\left[\pi^+\!\left(c(X)c(Y)c(Z)D_{\mathrm{dR}}^{-1}\right)\circ\pi^+\!\left(D_{\mathrm{dR}}^{-2m+2}\right)\right].
\]
\end{theorem}
The equality with the full residue follows from the vanishing of its
interior density.

\subsection{Boundary residue of \texorpdfstring{$\pi^+\!\bigl(c(X)c(Y)c(Z)D_{\mathrm{dR}}^{-2}\bigr)\circ\pi^+\!\bigl(D_{\mathrm{dR}}^{-2m+3}\bigr)$}{the de Rham--Hodge product with insertion c(X)c(Y)c(Z) and powers -2 and -2m+3}}\label{hd:secondsub}

The principal symbols of the two factors are now
\(c(X)c(Y)c(Z)/|\xi|^2\) and \(ic(\xi)(|\xi|^2)^{1-m}\).

\paragraph{Case \(a_1\).}
Take \(r=-2,\ \ell=-2m+3,\ j=k=0,\ |\alpha|=1\). The coefficient is \(-1\).
\begin{align*}
\widetilde\Phi_1={}&-\int_{|\xi'|=1}\int_{\mathbb R}\sum_{a=1}^{2m}
 \operatorname{tr}_{\Lambda}\!\Bigl[
 \partial_{\xi_a}\pi^+_{\xi_n}\sigma_{-2}\!\left(c(X)c(Y)c(Z)D_{\mathrm{dR}}^{-2}\right)\\
&\qquad{}\cdot \partial_{x_a}\partial_{\xi_n}\sigma_{-2m+3}\!\left(D_{\mathrm{dR}}^{-2m+3}\right)\Bigr](x_0,\xi)
 \,d\xi_n\,dS(\xi')\,d\operatorname{vol}_{\partial M}.
\end{align*}
At \(x_0\) both \(\partial_{x_a}c(\xi)\) and
\(\partial_{x_a}|\xi|^2\) vanish for every \(1\le a\le 2m\) and
every \(\xi_n\). Consequently,
\[
 \partial_{x_a}\partial_{\xi_n}\!
 \left[ic(\xi)(|\xi|^2)^{1-m}\right](x_0,\xi)=0,
\]
and the right factor in each summand vanishes.
The normal integral vanishes at each \(\xi'\), hence
\(\widetilde\Phi_1=0\).

\paragraph{Case \(a_2\).}
Take \(r=-2,\ \ell=-2m+3,\ j=1,\ k=|\alpha|=0\). The coefficient is \(-\frac12\).
\begin{align*}
\widetilde\Phi_2={}&-\frac12\int_{|\xi'|=1}\int_{\mathbb R}
 \operatorname{tr}_{\Lambda}\!\Bigl[
 \partial_{x_n}\pi^+_{\xi_n}\sigma_{-2}\!\left(c(X)c(Y)c(Z)D_{\mathrm{dR}}^{-2}\right)\\
&\qquad{}\cdot \partial_{\xi_n}^{2}\sigma_{-2m+3}\!\left(D_{\mathrm{dR}}^{-2m+3}\right)\Bigr](x_0,\xi)
 \,d\xi_n\,dS(\xi')\,d\operatorname{vol}_{\partial M}.
\end{align*}
The left projection throughout the collar is
\[
 -\frac{i c(X)c(Y)c(Z)}
 {2|\xi'|_g(\xi_n-i|\xi'|_g)}.
\]
Here both occurrences of \(|\xi'|_g\) contribute when taking the
normal derivative. Applying \eqref{hd:trace} after differentiating
the projected symbol gives the following expression. At fixed \(\xi'\) the trace before multiplication by
\(-1/2\) is
\begin{align*}
 &2^{2m+1}\sum_{j=1}^{n}\partial_{\xi_n}^{2}
       \left[\frac{\xi_j}{(1+\xi_n^2)^{m-1}}\right]\!
 \bigg\{\frac{\partial_{x_n}[g(X,Y)Z_j-g(X,Z)Y_j+g(Y,Z)X_j]}{2(\xi_n-i)}\\
 &\hspace{14mm}{}+h'(0)\left[-\frac1{4(\xi_n-i)}
           +\frac{i}{4(\xi_n-i)^2}\right]
          [g(X,Y)Z_j-g(X,Z)Y_j+g(Y,Z)X_j]\bigg\}.
\end{align*}
After multiplication by \(-1/2\), integrate twice by parts in \(\xi_n\)
to transfer both derivatives onto the pole factors. Applying
\eqref{eq:normal-moment} gives
\begin{align*}
 &\frac{6\pi(2m-2)!}{(m+1)!(m-2)!}\bigg\{
 \partial_{x_n}[g(X,Y)Z_n-g(X,Z)Y_n+g(Y,Z)X_n]\\
 &\qquad{}-\frac{5m h'(0)}{2(m+2)}[g(X,Y)Z_n-g(X,Z)Y_n+g(Y,Z)X_n]\\
 &+\frac{i(2m-1)}3\sum_{a=1}^{2m}\xi_a\,\partial_{x_n}[g(X,Y)Z_a-g(X,Z)Y_a+g(Y,Z)X_a]\\
 &-\frac{i(2m-1)(2m+1)h'(0)}{3(m+2)}
 \sum_{a=1}^{2m}\xi_a[g(X,Y)Z_a-g(X,Z)Y_a+g(Y,Z)X_a]\bigg\}.
\end{align*}
Only the first line survives angular integration, giving
\begin{align*}
 \widetilde\Phi_2={}&\frac{6\pi(2m-2)!\operatorname{Vol}(S^{2m-1})}{(m+1)!(m-2)!}\bigg\{
 \partial_{x_n}[g(X,Y)Z_n-g(X,Z)Y_n+g(Y,Z)X_n]\\
 &\qquad{}-\frac{5m h'(0)}{2(m+2)}[g(X,Y)Z_n-g(X,Z)Y_n+g(Y,Z)X_n]\bigg\}\,d\operatorname{vol}_{\partial M}.
\end{align*}

\paragraph{Case \(a_3\).}
Take \(r=-2,\ \ell=-2m+3,\ j=|\alpha|=0,\ k=1\). The coefficient is \(-\frac12\).
\begin{align*}
\widetilde\Phi_3={}&-\frac12\int_{|\xi'|=1}\int_{\mathbb R}
 \operatorname{tr}_{\Lambda}\!\Bigl[
 \partial_{\xi_n}\pi^+_{\xi_n}\sigma_{-2}\!\left(c(X)c(Y)c(Z)D_{\mathrm{dR}}^{-2}\right)\\
&\qquad{}\cdot \partial_{\xi_n}\partial_{x_n}\sigma_{-2m+3}\!\left(D_{\mathrm{dR}}^{-2m+3}\right)\Bigr](x_0,\xi)
 \,d\xi_n\,dS(\xi')\,d\operatorname{vol}_{\partial M}.
\end{align*}
The normal coordinate derivative of the right principal
symbol, taken before restriction, becomes
\[
 ih'(0)\left[\frac{c(\xi')}{2(1+\xi_n^2)^{m-1}}
 -\frac{(m-1)[c(\xi')+\xi_n c(dx_n)]}{(1+\xi_n^2)^m}\right].
\]
Its \(\xi_n\)-derivative is paired with the derivative of the left
projection, \(i c(X)c(Y)c(Z)/[2(\xi_n-i)^2]\).
Multiplying the trace by \(-1/2\) and integrating in \(\xi_n\) gives
\begin{align*}
 &\frac{6\pi(2m-2)!}{(m+1)!(m-2)!}\,h'(0)\bigg\{
 \frac{m(2m-1)}{2(m+2)}[g(X,Y)Z_n-g(X,Z)Y_n+g(Y,Z)X_n]\\
 &\qquad{}+\frac{i(2m-1)(m^2-1)}{3(m+2)}
 \sum_{a=1}^{2m}\xi_a[g(X,Y)Z_a-g(X,Z)Y_a+g(Y,Z)X_a]\bigg\}.
\end{align*}
Averaging over the tangential sphere therefore yields
\begin{align*}
 \widetilde\Phi_3={}&\frac{6\pi(2m-2)!\operatorname{Vol}(S^{2m-1})}{(m+1)!(m-2)!}\,\frac{m(2m-1)}{2(m+2)}\,h'(0)\\
 &\quad{}\times [g(X,Y)Z_n-g(X,Z)Y_n+g(Y,Z)X_n]\,d\operatorname{vol}_{\partial M}.
\end{align*}

\paragraph{Case \(b\).}
Take \(r=-3,\ \ell=-2m+3,\ j=k=|\alpha|=0\). The coefficient is \(-i\).
\begin{align*}
\widetilde\Phi_4={}&-i\int_{|\xi'|=1}\int_{\mathbb R}
 \operatorname{tr}_{\Lambda}\!\Bigl[
 \pi^+_{\xi_n}\sigma_{-3}\!\left(c(X)c(Y)c(Z)D_{\mathrm{dR}}^{-2}\right)\\
&\qquad{}\cdot \partial_{\xi_n}\sigma_{-2m+3}\!\left(D_{\mathrm{dR}}^{-2m+3}\right)\Bigr](x_0,\xi)
 \,d\xi_n\,dS(\xi')\,d\operatorname{vol}_{\partial M}.
\end{align*}
After the terms of zero trace
from \eqref{hd:square} are removed, its projected inverse-square factor
is
\begin{align*}
 -ih'(0)\bigg\{&
 -\frac12 c(\xi')c(dx_n)
       \left[-\frac1{4(\xi_n-i)^2}-\frac i{4(\xi_n-i)}\right]\\
 &-\frac{im}{4(\xi_n-i)^2}-\frac1{4(\xi_n-i)^3}
                       -\frac{i}{8(\xi_n-i)^2}\bigg\}.
\end{align*}
Multiplication on the left by \(c(X)c(Y)c(Z)\) restores the insertion.
Taking the trace against the normal derivative of the right symbol
and integrating with the factor \(-i\) gives
\begin{align*}
 &\frac{6\pi(2m-2)!}{(m+1)!(m-2)!}\,h'(0)\bigg\{
 \frac{m(m+7)}{2(m+2)}[g(X,Y)Z_n-g(X,Z)Y_n+g(Y,Z)X_n]\\
 &\qquad{}+\frac{i(2m-1)}3\left[-\frac1{m-1}+m+\frac{2m+1}{m+2}\right]\\
 &\hspace{29mm}{}\times\sum_{a=1}^{2m}\xi_a[g(X,Y)Z_a-g(X,Z)Y_a+g(Y,Z)X_a]\bigg\}.
\end{align*}
The term containing a single tangential coordinate has zero spherical mean, so
\begin{align*}
 \widetilde\Phi_4={}&\frac{6\pi(2m-2)!\operatorname{Vol}(S^{2m-1})}{(m+1)!(m-2)!}\,\frac{m(m+7)}{2(m+2)}\,h'(0)\\
 &\quad{}\times [g(X,Y)Z_n-g(X,Z)Y_n+g(Y,Z)X_n]\,d\operatorname{vol}_{\partial M}.
\end{align*}

\paragraph{Case \(c\).}
Take \(r=-2,\ \ell=-2m+2,\ j=k=|\alpha|=0\). The coefficient is \(-i\).
\begin{align*}
\widetilde\Phi_5={}&-i\int_{|\xi'|=1}\int_{\mathbb R}
 \operatorname{tr}_{\Lambda}\!\Bigl[
 \pi^+_{\xi_n}\sigma_{-2}\!\left(c(X)c(Y)c(Z)D_{\mathrm{dR}}^{-2}\right)\\
&\qquad{}\cdot \partial_{\xi_n}\sigma_{-2m+2}\!\left(D_{\mathrm{dR}}^{-2m+3}\right)\Bigr](x_0,\xi)
 \,d\xi_n\,dS(\xi')\,d\operatorname{vol}_{\partial M}.
\end{align*}
The right subleading symbol is that of
\((D_{\mathrm{dR}}^2)^{1-m}D_{\mathrm{dR}}\). With the mixed hatted
terms suppressed under the trace, this symbol is
\begin{align*}
 &h'(0)\left[\frac{(m-1)(m-\tfrac12)\xi_n}{(1+\xi_n^2)^m}
           +\frac{m(m-1)\xi_n}{(1+\xi_n^2)^{m+1}}\right]c(\xi')\\
 &\quad{}+h'(0)\left[-\frac{m}{2(1+\xi_n^2)^{m-1}}
 +\frac{(m-1)(m\xi_n^2-\tfrac12)}{(1+\xi_n^2)^m}
 +\frac{m(m-1)\xi_n^2}{(1+\xi_n^2)^{m+1}}\right]c(dx_n).
\end{align*}
Integrating by parts in case \(c\) changes \(-i\) into \(+i\).
Multiplication by the derivative of the left projection therefore
supplies the scalar factor \(-1/[2(\xi_n-i)^2]\) in front of the trace
of the insertion times this symbol.
The normal moments of the two Clifford coefficients at fixed \(\xi'\) now give
\begin{align*}
 &-\frac{6\pi(2m-2)!}{(m+1)!(m-2)!}\,h'(0)\bigg\{
 \frac{m(4m+3)}{2(m+2)}[g(X,Y)Z_n-g(X,Z)Y_n+g(Y,Z)X_n]\\
 &\qquad{}+\frac{i(2m-1)(2m^2+2m-1)}{3(m+2)}
 \sum_{a=1}^{2m}\xi_a[g(X,Y)Z_a-g(X,Z)Y_a+g(Y,Z)X_a]\bigg\}.
\end{align*}
The angular integral leaves
\begin{align*}
 \widetilde\Phi_5={}&-\frac{6\pi(2m-2)!\operatorname{Vol}(S^{2m-1})}{(m+1)!(m-2)!}\,\frac{m(4m+3)}{2(m+2)}\,h'(0)\\
 &\quad{}\times [g(X,Y)Z_n-g(X,Z)Y_n+g(Y,Z)X_n]\,d\operatorname{vol}_{\partial M}.
\end{align*}

The five cases give the following density, with the normal derivative
again supplied by \(\widetilde\Phi_2\).

\begin{theorem}\label{hd:second}
The boundary density of the second de Rham--Hodge residue is
\begin{align*}
 &\frac{6\pi(2m-2)!\operatorname{Vol}(S^{2m-1})}{(m+1)!(m-2)!}
 \bigg\{\left.\partial_{x_n}[g(X,Y)Z_n-g(X,Z)Y_n+g(Y,Z)X_n]\right|_{\partial M}\\
 &\hspace{29mm}{}-\frac m2 h'(0)
 [g(X,Y)Z_n-g(X,Z)Y_n+g(Y,Z)X_n]\bigg\}\,d\operatorname{vol}_{\partial M}.
\end{align*}
Its boundary integral equals
\[
 \widetilde{\operatorname{Wres}}\!\left[\pi^+\!\left(c(X)c(Y)c(Z)D_{\mathrm{dR}}^{-2}\right)\circ\pi^+\!\left(D_{\mathrm{dR}}^{-2m+3}\right)\right].
\]
\end{theorem}

\subsection{Hypersurface volume variation and the mean curvature term}\label{hd:meaning}

The coefficient of the boundary volume density in either theorem is
globally defined. The inward distance determines \(\partial_{x_n}\) and
the normal components throughout the collar, and hence determines the
normal derivative of the entire scalar pairing once \(X,Y,Z\) are fixed.
For comparison with the integrated mean curvature, recall that
\(K=\sum_{a=1}^{2m}g(\nabla^L_{e_a}\partial_{x_n},e_a)=-m h'(0)\).
Under the collar identification, the volume density of the hypersurface
at distance \(x_n\) is \(h(x_n)^{-m}d\operatorname{vol}_{\partial M}\).
The square root of its ratio to the boundary volume density is
\(h(x_n)^{-m/2}\). The common expression in braces in the two theorems
therefore has the geometric form
\[
 \left.\partial_{x_n}\!\left[
 h(x_n)^{-m/2}
 \big(g(X,Y)Z_n-g(X,Z)Y_n+g(Y,Z)X_n\big)\right]\right|_{\partial M}.
\]
Thus variation of the hypersurface volume accounts for the coefficient
\(K/2\) of the undifferentiated pairing.
The two local densities are proportional even when the three fields vary.
In particular,
\begin{align*}
 &\widetilde{\operatorname{Wres}}\!\left[\pi^+\!\left(AD_{\mathrm{dR}}^{-1}\right)\circ\pi^+\!\left(D_{\mathrm{dR}}^{-2m+2}\right)\right]\\
 &\quad=-\frac{2m-1}{3}\,\widetilde{\operatorname{Wres}}\!\left[\pi^+\!\left(AD_{\mathrm{dR}}^{-2}\right)\circ\pi^+\!\left(D_{\mathrm{dR}}^{-2m+3}\right)\right].
\end{align*}
\begin{corollary}\label{hd:curvature}
Suppose that \(g(X,Y)Z_n-g(X,Z)Y_n+g(Y,Z)X_n=C_0\) on a collar, for a
constant \(C_0\). Then
\begin{align*}
 &\widetilde{\operatorname{Wres}}\!\left[\pi^+\!\left(c(X)c(Y)c(Z)D_{\mathrm{dR}}^{-1}\right)\circ\pi^+\!\left(D_{\mathrm{dR}}^{-2m+2}\right)\right]\\
 &\quad=-\frac{\pi(2m-1)!\operatorname{Vol}(S^{2m-1})}{(m+1)!(m-2)!}\,
 C_0\!\int_{\partial M}K\,d\operatorname{vol}_{\partial M},\\
 &\widetilde{\operatorname{Wres}}\!\left[\pi^+\!\left(c(X)c(Y)c(Z)D_{\mathrm{dR}}^{-2}\right)\circ\pi^+\!\left(D_{\mathrm{dR}}^{-2m+3}\right)\right]\\
 &\quad=\frac{3\pi(2m-2)!\operatorname{Vol}(S^{2m-1})}{(m+1)!(m-2)!}\,
 C_0\!\int_{\partial M}K\,d\operatorname{vol}_{\partial M}.
\end{align*}
\end{corollary}
Indeed, the normal derivative of the scalar pairing is then zero.
The choice \(X=Y=Z=\partial_{x_n}\) near the boundary gives \(C_0=1\)
on every boundary component. Each residue therefore recovers the
integrated extrinsic curvature that appears in the Gibbons--Hawking--York
boundary term \cite{GH}, with the indicated inward-normal convention and numerical
normalization. For general fields the additional normal derivative
records the first normal jet of their insertion.

\begin{corollary}\label{hd:tangent}
If \(X,Y,Z\) are tangential on \(\partial M\), the curvature terms vanish
and the expression inside braces in each theorem reduces to
\[
 \left.[g(X,Y)\partial_{x_n}Z_n-g(X,Z)\partial_{x_n}Y_n
                  +g(Y,Z)\partial_{x_n}X_n]\right|_{\partial M}.
\]
If these fields are tangential on every nearby level set of the normal
distance, both densities vanish pointwise.
\end{corollary}

\paragraph{Half-density transport in the normal direction.}
For brevity put $C=g(X,Y)Z_n-g(X,Z)Y_n+g(Y,Z)X_n$ and $J(x_n)=h(x_n)^{-m}$. The bracket in the two formulas is
\[
 \mathcal L_\nu^{1/2}C:=\partial_{x_n}C+\frac K2C
   =\left.J^{-1/2}\partial_{x_n}(J^{1/2}C)\right|_{x_n=0}.
\]
Here the notation refers to normal transport of a scalar amplitude multiplied by the square root of the slice-volume density. The factor $K/2$ is forced by that square root. In particular, the condition $\mathcal L_\nu^{1/2}C=0$ means that the weighted amplitude is instantaneously constant in the normal direction; it need not mean that $C$ itself is constant. This gives a precise geometric explanation of the derivative and curvature appearing together.

\paragraph{Local mean curvature and first variation of area.}
The curvature measurement can be localized as well as integrated. Choose $X=f\nu$, $Y=Z=\nu$, with $f$ normally constant near the boundary. Then $C=f$ and $\partial_{x_n}C=0$. The two residues become their respective constants in Corollary~\ref{hd:curvature}, with $C_0\int K$ replaced by $\int fK$. Arbitrary smooth $f$ therefore recovers $K$ pointwise in the distributional, and hence smooth, sense. For the inward normal deformation of a closed boundary hypersurface with speed $f$, our sign convention gives
\[
 \left.\frac{d}{dt}\operatorname{Area}(\Sigma_t)\right|_{t=0}
     =\int_{\partial M}fK\,d\operatorname{vol}_{\partial M}.
\]
Thus normally constant normal probes measure the first variation of hypersurface area. A vanishing unweighted integral alone does not imply minimality; vanishing for every localized probe does imply $K=0$. In the present warped collar this also gives $h'(0)=0$, so the second fundamental form vanishes at the boundary.

\paragraph{Gravitational meaning and independence of spin structure.}
The same integral of extrinsic curvature occurs in the gravitational boundary term \cite{York,GH}. Its occurrence in spectral constructions has also been studied through boundary spectral actions \cite{CCB,CCBoundary}. Here it arises from an ordinary trace on all differential forms and therefore does not require a spin structure. The resulting identification is of a geometric density with a known explicit coefficient; it does not fix Newton's constant or establish the relative normalization of bulk and boundary terms in a gravitational variational problem. Indeed, the interior residue under consideration vanishes. This distinguishes the present boundary measurement from a derivation of the full Einstein--Hilbert action.

\paragraph{A torsion-free comparison and the limits of detection.}
The vanishing weighted sum is a useful comparison with Sections~\ref{sec:torsion} and~\ref{sec:connes}: the Levi--Civita de Rham operator supplies only the geometric expression, and its two coefficients are proportional. After division by fiber dimensions, the geometric coefficients agree with the spinor calculation, as follows directly from the displayed formulas. The full exterior module is therefore a natural reference for the metric contribution. Nevertheless, neither this proportionality nor the weighted cancellation is a topological index statement: an ordinary trace is used, and the surviving quantities vary with the metric and the first normal jets of the probes. They describe extrinsic local geometry rather than an Euler characteristic.

\end{document}